\documentclass[11pt]{article}

\usepackage[T1]{fontenc}

\usepackage{amssymb}
\usepackage{amsmath}

\usepackage{array}
\usepackage{hhline}
\usepackage{epsfig}
\usepackage{psfrag}
\usepackage{afterpage}
\usepackage{theorem}
\usepackage{stmaryrd}
\usepackage{mathdots}
\usepackage{color}
\usepackage{framed}

\makeatletter
\def\udots{\mathinner{\mkern1mu\raise\p@
\vbox{\kern7\p@\hbox{.}}\mkern2mu
\raise4\p@\hbox{.}\mkern2mu\raise7\p@\hbox{.}\mkern1mu}}
\makeatother

\usepackage{eufrak}
\newcommand{\goth}[1]{\EuFrak{#1}}
\newlength{\listeespacee}
\newlength{\tableautasse}
\newcommand{\ee}{{\goth{e}}}
\newcommand{\sss}{{\goth{s}}}

\newcommand{\HH}{\mbox{$\mathbb H$}}
\newcommand{\sign}{\operatorname{sign}}
\newcommand{\SO}{\operatorname{SO}}
\newcommand{\Sp}{\operatorname{Sp}}
\newcommand{\mat}{\mbox{\rm M}}

\newcommand{\T}{{\rm T}}

\newcommand{\dd}{\,\text{\rm d}}
\newcommand{\SSS}{{\cal S}}

\newcommand{\MM}{{\cal M}}

\newcommand{\KK}{{\cal K}}

\newcommand{\TT}{{\cal T}}

\newcommand{\GG}{{\cal G}}

\newcommand{\HHH}{{\cal H}}
\newcommand{\OO}{{\cal O}}

\newcommand{\R}{\mathbb R}
\newcommand{\N}{\mathbb N}

\newcommand{\C}{\mathbb C}
\newcommand{\K}{\mathbb K}

\newcommand{\Rad}{\operatorname{Rad}}
\newcommand{\codim}{\operatorname{codim}}
\newcommand{\Id}{\operatorname{Id}}
\newcommand{\Ric}{\operatorname{Ric}}
\newcommand{\ric}{\operatorname{ric}}

\newcommand{\Hbar}{\overline{H}}
\newcommand{\Zbar}{\overline{Z}}
\newcommand{\zbar}{\overline{z}}
\newcommand{\abar}{\overline{a}}

\newcommand{\ubar}{\overline{u}}
\newcommand{\vbar}{\overline{v}}
\newcommand{\wbar}{\overline{w}}

\newcommand{\tr}{\operatorname{tr}}

\newcommand{\End}{\operatorname{End}}

\newcommand{\Mat}{\mbox{\rm Mat}}
\newcommand{\vect}{\operatorname{vect}}

\newcommand{\im}{\operatorname{Im}}
\newcommand{\diag}{\operatorname{diag}}
\newcommand{\Span}{\operatorname{span}}

\newcommand{\sym}{+}
\newcommand{\asym}{-}

\newlength{\aux}

\newlength{\arraycolsepsauvegardegenerale}
\newlength{\arraycolsepsauvegardetemporaire}

\newlength{\tabcolsepsauvegardegenerale}
\newlength{\tabcolsepsauvegardetemporaire}

{\setlength{\tabcolsepsauvegardetemporaire}{\tabcolsep}\setlength{\tabcolsep}{0.5\tabcolsepsauvegardegenerale}\begin{tabular}}
{\end{tabular}\setlength{\tabcolsep}{\tabcolsepsauvegardetemporaire}}

{\setlength{\tabcolsepsauvegardetemporaire}{\tabcolsep}\setlength{\tabcolsep}{0.33\tabcolsepsauvegardegenerale}\begin{tabular}}
{\end{tabular}\setlength{\tabcolsep}{\tabcolsepsauvegardetemporaire}}

\newenvironment{narray}%
{\setlength{\arraycolsepsauvegardetemporaire}{\arraycolsep}\setlength{\arraycolsep}{0.5\arraycolsepsauvegardegenerale}\begin{array}}
{\end{array}\setlength{\arraycolsep}{\arraycolsepsauvegardetemporaire}}

\newenvironment{nnarray}%
{\setlength{\arraycolsepsauvegardetemporaire}{\arraycolsep}\setlength{\arraycolsep}{0cm}\begin{array}}
{\end{array}\setlength{\arraycolsep}{\arraycolsepsauvegardetemporaire}}

{\begin{list}{--}{\setlength{\leftmargin}{0cm}\setlength{\itemindent}{1.05cm}\setlength{\labelwidth}{1cm}\setlength{\labelsep}{0.5em}\setlength{\listparindent}{\parindent}}}%
{\end{list}}
{\begin{list}{}{\setlength{\leftmargin}{0cm}\setlength{\itemindent}{0cm}\setlength{\labelwidth}{1.5cm}\setlength{\labelsep}{0.6cm}\setlength{\listparindent}{\parindent}}}%
{\end{list}}

\theoremstyle{change}
\newtheorem{enonce}{}[section]

\newtheorem{de}[enonce]{Definition}
\newtheorem{prop}[enonce]{Proposition}
\newtheorem{prop-de}[enonce]{Proposition/Definition}
\newtheorem{de-prop}[enonce]{Definition/Proposition}
\newtheorem{lem}[enonce]{Lemma}
\newtheorem{te}[enonce]{Theorem}
\newtheorem{cor}[enonce]{Corollary}

\theorembodyfont{\rmfamily}
\newtheorem{rem-notation}[enonce]{Remark/Notation}
\newtheorem{lem-notation}[enonce]{Lemma/Notation}
\newtheorem{def-notation}[enonce]{Definition/Notation}
\newtheorem{notation}[enonce]{Notation}
\newtheorem{rem}[enonce]{Remark}
\newtheorem{rem-de}[enonce]{Remark/Definition}
\newtheorem{importantrem}[enonce]{Important Remark}

\newtheorem{terminology}[enonce]{Terminology}
\newtheorem{terminology-notation}[enonce]{Terminology/Notation}

\newtheorem{consequence-notation}[enonce]{Consequence/Notation}
\newtheorem{reminder}[enonce]{Reminder}
\newtheorem{reminder-rem}[enonce]{Reminder/Remark}

\newcounter{claimcounter}
\psfrag{1}{{\small\bf (1)}}
\psfrag{1c}{{\small\bf (1$^\C$)}}
\psfrag{2}{{\small\bf (2)}}
\psfrag{2'}{{\small\bf (2')}}
\psfrag{2c}{{\small\bf (2$^\C$)}}
\psfrag{3}{{\small\bf (3)}}
\psfrag{3'}{{\small\bf (3')}}
\psfrag{3c}{{\small\bf (3$^\C$)}}

\author{Charles Boubel}

\begin{document}
\thispagestyle{empty}
\begin{center}
{\bf THE ALGEBRA OF PARALLEL ENDOMORPHISMS OF A
PSEUDO-RIEMANNIAN METRIC: SEMI-SIMPLE PART\bigskip\\}
{\sc Charles Boubel}\footnote{Institut de Recherche Math\'ematique Avanc\'ee, UMR 7501 -- 
Universit\'e de Strasbourg et CNRS, 7 rue Ren\'e Descartes, 67084 STRASBOURG CEDEX, FRANCE}\medskip\\November 29th.\@ 2013, revised April 16th. 2015
\end{center}

\begin{center}\setlength{\aux}{\textwidth}\addtolength{\aux}{-1.7cm}
\parbox{\aux}{\small {\bf Abstract.} On a (pseudo\nobreakdash-)Riemannian manifold $(\MM,g)$, some fields of endomorphisms {\em i.e.\@} sections of $\End(T\MM)$ may be parallel for $g$. They form an associative algebra $\goth e$, which is also the commutant of the holonomy group of $g$. As any associative algebra, $\goth e$ is the sum of its radical and of a semi-simple algebra $\goth s$. Here we study $\goth s$: it may be of eight different types, including the generic type $\goth s=\R\Id$, and the K\"ahler and hyperk\"ahler types $\goth s\simeq\C$ and $\goth s\simeq{\mathbb H}$. This is a result on real, semi-simple algebras with involution. For each type, the corresponding set of germs of metrics is non-empty; we parametrise it. We give the constraints imposed to the Ricci curvature by parallel endomorphism fields.
\medskip

\medskip

\noindent{\bf M.S.C.\@ 2010:} 53B30, 53C29, 16K20, 16W10 secondary 53B35, 53C10, 53C12, 15A21.}
\end{center}

We classify here the germs of (pseudo\nobreakdash-)Riemannian metric, after the semi-simple part of their algebra of parallel endomorphism fields. Our motivation is the following.\\

\noindent{\bf Motivation.} A K\"ahler metric $g$ on some manifold $\MM$ may be defined as a Riemannian metric admitting an almost complex structure $J$ which is parallel: $DJ=0$ with $D$ the Levi-Civita connection of $g$. A natural question is to ask whether other fields of endomorphisms, {\em i.e.\@} sections of $\End(T\MM)$, may be parallel for a Riemannian metric. The answer is nearly immediate. First, one restricts the study to metrics that do not split into a non trivial Riemannian product, called here ``indecomposable''. Otherwise, any parallel endomorphism field is the direct sum of parallel such fields on each factor (considering as a unique factor the possible flat factor). Then a brief reasoning ensures that only three cases occur: $g$ may be generic {\em i.e.\@} admit only the homotheties as parallel endomorphisms, be K\"ahler, or be hyperk\"ahler {\em i.e.\@} admit two (hence three) anticommuting parallel complex structures. The brevity of this list is due to a simple fact: the action of the holonomy group $H$ of an indecomposable Riemannian metric is irreducible {\em i.e.\@} does not stabilise any proper subspace. In particular, this compels any parallel endomorphism field to be of the form $\lambda\Id+\mu J$ with $J$ some parallel, skew adjoint almost complex structure. Now, such irreducibility fails in general for an indecomposable pseudo-Riemannian metric, so that a miscellany of other parallel endomorphism fields may appear. This gives rise to the question tackled here:
\begin{center}
Which (algebra of) parallel endomorphism fields may a pseudo-Riemannian metric admit ?
\end{center}

The interest of this question lies also in the following. When studying the holonomy of indecomposable pseudo-Riemannian metrics, the irreducible case may be exhaustively treated: the full list of possible groups, together with the corresponding spaces of germs of metrics (and possibly compact examples) may be provided. After a long story that we do not recall here, this has been done, even for germs of arbitrary torsion free affine connections, see {\em e.g.\@} the surveys \cite{bryant1996,schwachhoefer}. Yet, in general, the representation of $H$ may be non-semi-simple and such an exhaustive answer is out of reach, except perhaps in very low dimension, see {\em e.g.\@} the already long list of possible groups in dimension four in \cite{BB-ikemakhen1997,ghanam-thompson2001}.
Thus, intermediate questions are needed: not aiming at the full classification, but still significant; see {\em e.g.} \cite{galaev-leistner2010} for a survey of such works. Investigating the commutant $\End(T_m\MM)^H$ of $H$ at some point $m$ of $\MM$, instead of $H$ itself ---~that is to say studying the algebra of parallel endomorphisms~--- is such a question.
One may also notice that determining all the parallel tensors, not only the endomorphisms, would mean determining the algebraic closure of the holonomy group $H$. So this work is a step towards this.

Now, as any associative algebra, $\End(T_m\MM)^H$ classically splits into a sum $\goth{s}\oplus\goth{n}$ with $\goth s$ a semi-simple subalgebra ---~in general not canonical, but its isomorphism class is~--- and $\goth n:=\Rad(\End(T_m\MM)^H)$ a nilpotent ideal, its radical. The study of $\goth s$ and $\goth n$ involve very different methods, and each of them is a work in itself. This article is devoted to $\goth s$; we deal with $\goth n$ in \cite{boubel13} and other future works. The interest of this article is that:

\begin{center} 
We deal with indecomposable metrics the holonomy group of which is never supposed to be irreducible or totally reducible.
\end{center}

As it is classical in holonomy problems, the question is twofold: {\bf (i)} Which algebras $\goth s$ are possible ? {\bf (ii)} By which sets of metrics are they produced ? We will handle both, restraining ourselves, for point {\bf (ii)}, to the first natural step {\em i.e.\@} to {\em germs} of metrics.\\

\noindent{\bf Contents and structure of the article.} Let $(\MM,g)$ be a (pseudo\nobreakdash-)Riemannian manifold of dimension $d$, $H$ its holonomy group, $H^0$ the neutral component of $H$ and $m\in\MM$.\medskip

In Part 1, we introduce the decomposition $\End(T_m\MM)^H=\goth s\oplus\goth n$ in \S\ref{s_et_n} and some natural objects associated with a {\em reducible} holonomy representation in \S\ref{considerations}, together with a simple but remarkable commutation property in $\End(T_m\MM)^H$, Proposition \ref{pseudocommutation}. In \S\ref{classification_s}, we give {\bf our main theorem}: $\goth s$ may be of eight different types, including the generic, K\"ahler and hyperk\"ahler types $\goth s=\R\Id$, $\goth s\simeq\C$ and $\goth s\simeq\HH$. See Theorem \ref{structure_s} p.\@ \pageref{structure_s} and Tables \ref{table_1} and \ref{table_2} for details. In the five non-Riemannian cases, the metric has necessarily a ``neutral'' signature $(\frac d2,\frac d2)$ and $\goth s$ contains a ``parak\"ahler'' structure $L$ {\em i.e.\@} a $g$-skew adjoint automorphism such that $T\MM=\ker(N-\Id)\oplus\ker(N+\Id)$. (Such a structure is also called a ``bi-Lagrangian'' structure, see Terminology \ref{bilagrangian}.) That is linear algebra: the classification of some semi-simple, $g$-self adjoint subalgebras of ${\goth gl}(\R^d)$, see Remark \ref{resultatdalgebrelineaire}. We also give two corollaries.\medskip

In Part 2, to show that each type given by Theorem \ref{structure_s} occurs, we do a little more: we parametrise the set of germs of metrics in each of them (explicitly, or {\em via} Cartan-K\"ahler theory). Here we adapt a classical proof the line of which is given by R.\@ Bryant in \cite{bryant1996} ---~in particular, in the Riemannian case, this provides an explicit writing of this proof.\medskip

In Part 3, we give the consequences of the existence of any type of parallel endomorphisms on the Ricci curvature. They are quite simple, hence very remarkable.\\

\noindent{\bf General setting and some general notation.} Here $\MM$ is a simply connected manifold of dimension $d$ and $g$ a Riemannian or pseudo-Riemannian metric on it, whose holonomy representation does not stabilise any nondegenerate subspace, that is to say does not split in an orthogonal sum of subrepresentations. In particular, $g$ does not split into a Riemannian product. We set $H\subset\SO^0(T_m\MM,g_{|m})$ the holonomy group of $g$ at $m$ and $\goth h$ its Lie algebra. As $\MM$ is supposed to be simply connected, we deal everywhere with $\goth h$, forgetting $H$. Let $\goth{e}$ be the algebra $\End(T_m\MM)^{\goth{h}}$  of the parallel endomorphisms of $g$ ---~to commute with $\goth h$ amounts to extend as a parallel field~---; it is isomorphic to some subalgebra of $\mat_d(\R)^{\goth{h}}$. Notice that $\ee$ is stable by $g$-adjunction, which we denote by $\sigma:a\mapsto a^\ast$. If $A$ is an algebra and $B\subset A$, we denote by $\langle B\rangle$, $(B)$, and $A^B$ the algebra, respectively the ideal, spanned by $B$, and the commutant of $B$ in $A$. When lower case letters: $x_i, y_i$ {\em etc.\@} stand for local coordinates, the corresponding upper case letters: $X_i, Y_i$  {\em etc.\@} stand for the corresponding coordinate vector fields. Viewing vector fields $X$ as derivations, we denote Lie derivatives ${\cal L}_Xu$ also by $X.u$.

The matrix $\diag(I_p,-I_q)\in\mathrm{M}_{p+q}(\R)$ is denoted by \addtolength{\arraycolsep}{-.5ex}$I_{p,q}$, $\mbox{\footnotesize$\left(\begin{array}{cc}0&-I_{p}\\I_{p}&0\end{array}\right)$}\in\mathrm{M}_{2p}(\R)$\addtolength{\arraycolsep}{.5ex} by $J_p$ and $\mbox{\footnotesize$\left(\begin{array}{cc}0&I_{p}\\I_{p}&0\end{array}\right)$}\in\mathrm{M}_{2p}(\R)$ by $L_p$. If $V$ is a vector space of even dimension $d$, we recall that an $L\in\End(V)$ is called {\em paracomplex} if $L^2=\Id$ with $\dim\ker(L-\Id)=\dim\ker(L+\Id)=\frac d2$.

Finally, take $A\in\Gamma(\End(T\MM))$, paracomplex {\em i.e.\@} such that $\dim\ker(N-\Id)=\dim\ker(N+\Id)=\frac d2$. If it is integrable {\em i.e.\@} if its matrix is constant in well-chosen local coordinates, we call it a ``paracomplex structure'', like a complex structure, as opposed to an almost complex one.\\

\noindent{\bf Acknoledgements.} I thank M.\@ Brion for a few crucial pieces of information and references in Algebra, P.\@ Baumann for his availability and for the references he indicated to me. I thank M.\@ Audin, P.\@ Mounoud and P.\@ Py for their comments on the writing of certain parts of the manuscript, and the referee for a few useful remarks.

\section{\mathversion{bold}The algebra $\goth e=\End(T\MM)^H$ and its semi-simple part $\goth s$\mathversion{normal}}\label{semisimple}

\subsection{The decomposition $\goth e=\goth{s}\oplus\goth{n}$ of $\goth e$ in a semi-simple part and its radical }\label{s_et_n}

First we need to recall some facts and set some notation. All the results invoked are classical for finite dimensional associative algebras; we state them for a unital real algebra $A$.

\begin{notation}
If $A$ is a subset of an algebra, $A^\ast\subset A$ denotes here the subset of its invertible elements. If $\sigma$ is an involutive anti morphism of $A$, then $A^\pm=\{U\in A;\sigma(A)=\pm A\}$ denotes the subspace of its self adjoint or skew adjoint elements.
\end{notation}

\begin{reminder}\label{n_nilpotent} An algebra $A$ is said to be {\em nilpotent} if $A^k$, the algebra spanned by the products of $k$ elements of $A$, is $\{0\}$ for some $k$. In particular, the elements of a nilpotent subalgebra of ${\rm M}_n(\R)$ are simultaneously strictly upper triangular in some well-chosen basis.
\end{reminder}

\begin{de} (See \cite{curtis-reiner} \S25 or \cite{jacobson_structureofrings}) The radical $\Rad A$ of $A$ is the intersection of its maximal ideals. It is a nilpotent ideal. Equivalently, it is the sum of the nilpotent ideals of $A$. The algebra $A$ is said to be simple if its only proper ideal is $\{0\}$, and semi-simple if its radical is $\{0\}$ ---~so a simple algebra is semi-simple, and $A/\Rad(A)$ is semi-simple.
\end{de}

The decomposition $\goth e=\goth{s}\oplus\goth{n}$ is provided by the following classical result. The last assertion is a refinement due to Taft \cite{taft1, taft2}. I thank P.\@ Baumann for this reference.

\begin{te}\label{wm} {\bf [Wedderburn -- Mal\mathversion{bold}$\check{\mathbf{c}}$\mathversion{normal}ev]} (see \cite{curtis-reiner} \S72) Let $A$ be a finite dimensional $\R$-algebra. Then there exists a semi-simple algebra $A_S$ in $A$ such that $A=A_S\oplus \Rad(A)$. If moreover $A$ is endowed with an involutive anti-morphism $\sigma$, then $A_S$ may be chosen $\sigma$-stable.
\end{te}

\begin{notation}We set $\goth n=\Rad \goth e$. Being the unique maximal nilpotent ideal of $\ee$, $\goth n$ is self adjoint {\em i.e.} stable by $g$-adjunction. We take $\goth s\simeq\goth{e}/\goth{n}$ some self adjoint semi-simple subalgebra of $\goth e$ provided by Theorem \ref{wm}.
\end{notation}

\subsection{Some natural objects associated with a reducible holonomy representation; a ``quasi-commutation'' property}\label{considerations}

\begin{rem-notation} We denote by $E_0=\cap_{W\in\goth{h}}\ker W$ the (possibly trivial) maximal subspace of $T_m\MM$ on which the holonomy group $H$ acts trivially. As $T_m\MM$ is $H$-orthogonally indecomposable, $E_0$ is totally isotropic. We set ${\goth n}_0=\{N\in\goth{e}\,;\,\im N\subset E_0\}$; as the actions of $H$ and $\goth e$ on $T_m\MM$ commute, ${\goth n}_0$ is an ideal of $\goth e$, moreover self adjoint. So, for any $x,y\in T_m\MM$, and any $N,N'\in\goth{n}_0$, $g(N'Nx,y)=g(Nx,N^{\prime\ast}y)\in g(E_0,E_0)=\{0\}$, so $N'N=0$ {\em i.e.\@} $\goth{n}_0^2=\{0\}$.
\end{rem-notation}

\begin{rem-notation}\label{produit_naturel} The algebra $\goth e$ is naturally endowed with the bilinear symmetric form $\langle U,V\rangle=\frac1d\tr(U^\ast V)$. By Reminder \ref{n_nilpotent}, $\goth n\subset\ker(\langle\,\cdot\,,\,\cdot\,\rangle)$. If moreover $\goth e$ admits some self adjoint complex structure $\underline J$, and denoting by $\goth{e}_{\underline J}$ the $\underline J$-complex algebra $\{U\in\goth e\,;\,U{\underline J}={\underline J}U\}$, then $\goth{e}_{\underline J}$ is endowed with the complex form $\langle U,V\rangle_{\underline J}=\frac1d(\tr(U^\ast V)-{\mathrm i}\tr(U^\ast {\underline J}V))$.
\end{rem-notation} 

The following proposition is the key of most steps of the classification \ref{structure_s}. As it is also worth to be noticed by itself, we state it apart, here.

\begin{prop}\label{pseudocommutation} Let $U,V$ be in $\goth{e}$ and $m$ be any point of $\MM$. If $U$ is self adjoint, then for any $x,y\in T_m\MM$, $R(x,y)(UV-VU)=0$. Consequently, $UV-VU\in\goth n_0$. In particular, in case $E_0=\cap_{W\in\goth{h}}\ker W$ is reduced to $\{0\}$, all self adjoint elements of $\goth{e}$ are central in $\goth{e}$.
\end{prop}

Proposition \ref{pseudocommutation} rests on the following remark.

\begin{reminder-rem}\label{g_U}Classically, the Bianchi identity implies that, at all point $m\in\MM$:
\begin{equation}\label{swap}
\forall x,y,z,t\in T_m\MM,\ g(R(x,y)z,t)=g(R(z,t)x,y).
\end{equation}
This holds also if we replace $g$ by any bilinear form parallel with respect to the Levi Civita connection of $g$, degenerate or not. The proof does not need nondegeneracy, see {\em e.g.\@} Lemma 9.3 in \cite{milnor_morse}. A consequence is that, if $U$ is a parallel self adjoint endomorphism:
$$\forall x,y,z\in T_m\MM,\ R(Ux,y)z=R(x,Uy)z=R(x,y)Uz.$$
The first equality is classical. For the second one, take $t$ any fourth vector and denote by $g_U$ the bilinear form $g(\,\cdot\,,U\,\cdot\,)$, which is parallel, as $U$ is, and symmetric, as $U^\ast=U$. Then:
\begin{align*}
g(R(x,Uy)z,t)&=g(R(z,t)x,Uy)\qquad\text{applying (\ref{swap}) to $g$,}\\
&=g_U(R(z,t)x,y)\\
&=g_U(R(x,y)z,t)\qquad\text{applying (\ref{swap}) to $g_U$,}\\
&=g(R(x,y)Uz,t)\qquad\text{as $U^\ast=U$, being parallel,
}\\&\qquad\qquad\qquad\qquad\qquad\text{commutes with $R(x,y)$.}\tag*{$\Box$}
\end{align*}
\end{reminder-rem}

\noindent{\bf Proof of Proposition \ref{pseudocommutation}.} Take $U,V\in\goth{e}$ with $U^\ast=U$ and $x,y,z,t\in T_m\MM$. The bilinear form $g_U:=g(\,\cdot\,,U\,\cdot\,)$ is parallel, as $U$ is.
\begin{align*}
&g(R(x,y)z,VUt)\\
\ \ \ =&g(R(x,y)V^\ast z,Ut)\qquad\text{as, $V^\ast$, parallel, commutes with $R(x,y)$,}\\
=&g(R(x,Uy)V^\ast z,t)\qquad\text{by Remark \ref{g_U}, applied to $U$,}\\
=&g(R(x,Uy)z,Vt)\qquad\text{as, $V^\ast$ commutes with $R(x,y)$,}\\
=&g(R(x,y)z,UVt)\qquad\text{again by Remark \ref{g_U}, so the result.}\tag*{$\Box$}
\end{align*}

\subsection{The eight possible forms of $\goth s$}\label{classification_s}

The types given by Theorem \ref{structure_s} are known, but not in full generality for type {\bf (3')} {\em i.e.\@} with the corresponding set of germs of metrics clearly stated, and except {\bf (3$^\C$)} which I never encountered explicitly. So Theorem \ref{structure_s} closes the list, {\em may the action of $H$ be totally reducible or not}. The proof rests on the classical Wedderburn-Artin and Skolem Noether theorems, and then is elementary. Remark \ref{types_pour_h} below gives the generic holonomy group corresponding to each case of the theorem.

\begin{te}\label{structure_s}
The algebra $\sss$ is of one of the following types, where $\underline J$, $J$,  and $L$ denote respectively self adjoint complex structures and skew adjoint complex and paracomplex structures. Each case is precisely described in Tables \ref{table_1} p.\@ \pageref{table_1} and \ref{table_2} p.\@ \pageref{table_2}, which are part of the theorem.\medskip

\noindent{\bf  (1) generic}, $\sss=\vect(\Id)$.\medskip

\noindent{\bf  (1$^{\C}$) ``complex Riemannian''}, $\sss=\vect(\Id,\underline J)$. Here $d\geqslant4$ is even, $\sign(g)=(\frac{d}{2},\frac{d}{2})$, $(\MM,\linebreak[1]J,\linebreak[1]g(\cdot,\cdot)-\mbox{\rm i}g(\cdot,\underline J\cdot))$ is complex Riemannian for a unique complex structure in $\goth s$, up to sign.\medskip

\noindent{\bf  (2) (pseudo\nobreakdash-)K\"ahler}, $\sss=\vect(\Id,J)$. Here $d$ is even and $(\MM,J,g)$ is (pseudo\nobreakdash-)K\"ahler, for a unique complex structure in $\goth s$, up to sign.\medskip

\noindent{\bf  (2') parak\"ahler}, $\sss=\vect(\Id,L)$. Here $d$ is even, $\sign(g)=(\frac{d}{2},\frac{d}{2})$, $(\MM,L,g)$ is parak\"ahler, for a unique paracomplex structure in $\goth s$, up to sign.\medskip

\noindent{\bf  (2$^{\C}$) ``complex K\"ahler''}, $\sss=\vect(\Id,\underline J,L,J)$. Here $d\in4\N^\ast$, $\sign(g)=(\frac{d}{2},\frac{d}{2})$ and $(\MM,\underline J,J,L,g)$ is at once complex Riemannian, pseudo-K\"ahler, and parak\"ahler, on a unique way in $\goth s$, up to sign of each structure.\medskip

\noindent{\bf  (3) (pseudo\nobreakdash-)hyperk\"ahler}, $\sss=\vect(\Id,J_1,J_2,J_3)$. Here $d\in4\N^\ast$, $(\MM,J_1,J_2,g)$ is (pseudo\nobreakdash-)hyperk\"ahler, the set of K\"ahler structures in $\goth s$ being a 2-dimensional submanifold.\medskip

\noindent{\bf  (3') ``para-hyperk\"ahler''}, $\sss=\vect(\Id,J,L_1,L_2)$. Here $d\in4\N^\ast$, $\sign(g)=(\frac{d}{2},\frac{d}{2})$ and $(\MM,J,L_1,g)$ is at once pseudo-K\"ahler and parak\"ahler, the set of complex and of paracomplex structures in $\goth s$ being each  a 2-dimensional submanifold.\medskip

\noindent{\bf  (3$^{\C}$) ``complex hyperk\"ahler''}, $\sss=\vect(\Id,\underline J,J,L_1,L_2,\underline JJ,\underline JL_1,\underline JL_2)$  Here $d\in8\N^\ast$, $\sign(g)=(\frac{d}{2},\frac{d}{2})$ and $(\MM,g)$ is at once complex Riemannian (on a unique way up to sign in $\goth s$), and pseudo-K\"ahler and parak\"ahler. The sets of pseudo- or parak\"ahler structures are 2-dimensional $\underline J$-complex submanifolds of $\goth s$.\medskip

Each type is produced by a non-empty set of germs of metrics. On a dense open subset of them, for the $C^2$ topology, the holonomy group of the metric is the commutant ${\rm SO}^0(g)^{\goth s}$ of ${\goth s}$ in ${\rm SO}^0(g)$. Cases {\bf  (3)}, {\bf  (3')}, and {\bf  (3$^{\C})$} are Ricci-flat, see Theorem \ref{ricci} p.\@ \pageref{ricci}.
\end{te}

\begin{terminology}\label{bilagrangian} A pseudo-riemannian metric $g$ with a parak\"ahler structure $L$ amounts to the data of a symplectic form $\omega$, of two distinguished transverse Lagrangian distributions $E^\pm$ and of a torsion free connection $D$ with respect to which $\omega$, $E^+$ and $E^-$ are parallel: set $\omega:=g(\,\cdot\,,L\,\cdot\,)$, $E^\pm=\ker(L\pm{\rm Id})$ and take for $D$ the Levi-Civita connection of $g$. Such a structure is also known as a ``bi-Lagrangian'' structure; this terminology has been introduced by R.\@ Bryant in \cite{bryant2001}, \S5.2. I thank the referee for this remark.
\end{terminology}

\begin{rem}The fact that the set of germs of metrics in each case is non-empty is well-known, except perhaps for types {\bf (3')} and {\bf (3$^{\C}$)}. In all cases, \S\ref{realisation} gives their parametrisation.
\end{rem}

\begin{rem}\label{resultatdalgebrelineaire}In fact, we proved the following result in plain linear algebra. If $g$ is a (pseudo\nobreakdash-)Euclidean product on $\R^d$ and $A$ a semi-simple, $g$-self adjoint subalgebra of $\goth{gl}(\R^d)$, whose action on $\R^d$ is indecomposable (in an orthogonal sum), then $A$ is one of the eight algebras of Theorem  \ref{structure_s} or the algebra $A\simeq\HH\oplus\HH$ of Remark \ref{bianchi_indispensable}.
\end{rem}

\begin{notation}\label{notationV}If $G$ is a subgroup of GL$_d(\K)$, we denote here by ${\bf V}$ its standard representation in $\K^d$. We denote then by ${\bf V}^\ast:g\mapsto(\lambda\mapsto\lambda\circ g^{-1})$ its contragredient representation in $(\K^d)^\ast$ and, if $\K=\C$, by $\overline{\bf V}^\ast$ the complex conjugate of it.
\end{notation}

\begin{rem}\label{representation_double} In cases {\bf (2')}, {\bf (2$^\C$)}, {\bf (3')} and {\bf (3$^\C$)}, the existence of a paracomplex structure $L$ splits $\T\MM=\ker(L-\Id)\oplus\ker(L+\Id)=V\oplus V'$ into a sum of two totally isotropic factors, and the morphism $\flat$ given by the metric identifies $V'$ with $V^\ast$. Then, $H$ is isomorphic to a subgroup $[H]$ of GL$_{d/2}(\K)$, the holonomy representation being ${\bf V}\oplus{\bf V}^\ast$,  if $\K=\R$, or ${\bf V}\oplus\overline{\bf V}^\ast$,  if $\K=\C$, on $\ker(L-\Id)\oplus\ker(L+\Id)$. Matricially:
$$H:=\left\{\text{\small$\left(\begin{array}{cc}U&0\\0&^t\overset{\raisebox{-.3ex}{$\scriptscriptstyle($}\rule{.8em}{.2mm}\raisebox{-.3ex}{$\scriptscriptstyle)$}}{U}^{\raisebox{-.7ex}{$\scriptscriptstyle-1$}}\end{array}\right)$},U\in [H]\right\};$$
so if $\K=\R$, $H\subset{\rm SO}^0(\frac d2,\frac d2)$ and if $\K=\C$, $H\subset{\rm U}(\frac d2,\frac d2)$.
\end{rem}

\begin{rem}\label{types_pour_h}For $\goth s$ of each type, we sum up here: the possible signature(s) of $g$, the group in which $H$ (possibly identified with $[H]$, see Rem.\@ \ref{representation_double}) is included, and to which it is generically equal (proof in \S\ref{realisation})
, and the representation of $H$ or $[H]$. Notice that each time, this group is also the commutant of ${\goth s}$ in ${\rm SO}^0(g)$. See Notation \ref{notationV} for {\bf V}.\medskip

\noindent\hspace*{\fill}{\begin{tabular}{cccccccc}
{\bf (1)}&{\bf (1$^\C$)}&{\bf (2)}&{\bf (2')}&{\bf (2$^\C$)}&{\bf (3)}&{\bf (3')}&{\bf (3$^\C$)}\\\hline
$(p,q)$&$(p,p)$&$(2p,2q)$&$(p,p)$&$(2p,2p)$&$(4p,4q)$&$(2p,2p)$&$(4p,4p)$
\\\hline
SO$^0(p,q)$&SO$(p,\C)$&U$(p,q)$&GL$^0(p,\R)$&GL$(p,\C)$&Sp$(p,q)$&Sp$(2p,\R)$&
Sp$(2p,\C)$
\\\hline
${\bf V}$&${\bf V}$&${\bf V}$&${\bf V}\oplus{\bf V}^\ast$&${\bf V}\oplus\overline{\bf V}^\ast$&${\bf V}$&${\bf V}\oplus{\bf V}^\ast$&${\bf V}\oplus\overline{\bf V}^\ast$
\end{tabular}}\hspace*{\fill}
\end{rem}

\begin{rem}In Theorem \ref{structure_s}, the new cases with respect to the Riemannian framework occur only for metrics $g$ of signature $\bigl(\frac d2,\frac d2\bigr)$.
\end{rem}

\begin{rem}{\bf [Justification of the labels in Theorem \ref{structure_s}]} The generic holonomy groups corresponding to $\goth s$ of types {\bf (1$^\C$)}, {\bf (2$^\C$)} and {\bf (3$^\C$)} are complexification of those corresponding to $\goth s$ of respective types {\bf (1)}, ({\bf (2)} or {\bf (2')}), and ({\bf (3)} or {\bf (3')}). Besides, if you consider the different types in a comprehensive sense, type {\bf (2)} {\em e.g.\@} meaning only ``$H\subset{\rm U}(p,q)$'', and so on, you obtain the following inclusion diagram:\medskip

\noindent\hspace*{\fill}\mbox{\epsfig{file=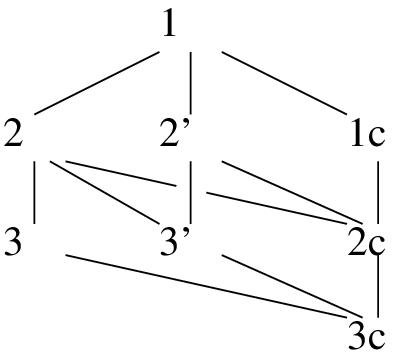,width=4cm}}\hspace*{\fill}\raisebox{1.5cm}{\parbox{5cm}{where the strokes denote the inclusion of the set of metrics below into the one above.}}\hspace*{\fill}\medskip

\noindent This justifies our notation. Another point of view is the following. Suppose that $g$ is a real analytic germ of metric at $m$. Then $\goth h$ is generated by $\{D^kR(u_1,\ldots,u_{k+2}),(u_i)_{i=1}^{k+2}\in T_m\MM\}$, the curvature tensors at $m$ and their covariant derivatives at all orders. So the complexification $g^\C$ of the germ $g$ has $\goth h\otimes\C$ as holonomy algebra. Thus here, if $g$ is ``of type {\bf (1)}'', respectively ({\bf (2)} or {\bf (2')}), or  ({\bf (3)} or {\bf (3')}), its complexification is ``of type {\bf (1$^\C$)}'', respectively {\bf (2$^\C$)} and {\bf (3$^\C$)}.
\end{rem}

\begin{table}[ht!]
$
\hspace*{\fill}\begin{array}{|c|cccc@{\hspace{-.ex}}c|}
\multicolumn{1}{c}{}&\multicolumn{1}{c}{\text{\bf (1)}}&\multicolumn{1}{c}{\text{\bf (1$^{\C}$)}}&\multicolumn{1}{c}{\text{\bf (2)}}&\multicolumn{1}{c}{\text{\bf (2$^{\prime}$)}}&\multicolumn{1}{c}{\text{\bf (2$^{\C}$)}}\\
\hline
\goth{s}=& \langle\ \rangle & \langle \underline J\rangle & \langle J\rangle & \langle L\rangle &\begin{nnarray}{c}{\langle \underline J,L,J\,|}\underline J\in Z(\goth{s}),\\LJ=JL=\underline J\rangle\end{nnarray}\\
\hline
\goth{s}\simeq & \R & \C & \C & \R\oplus\R & \C\oplus\C\\
\hline
\begin{nnarray}{c}
\text{through the}\\\text{morphism}\\\text{$\varphi$ given by:}\end{nnarray}&{\rm Id}\mapsto 1&\underline J\mapsto {\rm i}&J \mapsto {\rm i}&L\mapsto (1,-1)&\begin{nnarray}{c}(\underline J,L,J)\mapsto((\text{i},\text{i})\\(1,-1),(\mathrm{i},-\mathrm{i}))\end{nnarray}\\
\hline
\begin{nnarray}{c}\text{$\varphi$ conjugates}\\\text{adjunction}\\\text{in $\goth s$ to:}\end{nnarray}&\Id_\R&\Id_\C&z\mapsto\overline z&(a,b)\mapsto(b,a)&(a,b)\mapsto(b,a)\\
\hline
\varphi(\goth s^+)=&\R&\C&\R&\R.(1,1)&\C.(1,1)\\
\hline
\begin{nnarray}{c}\varphi(\goth s^+)\subset\varphi(\goth s)\\[-1ex]\text{\scriptsize[1]}\end{nnarray}&\R^{1,0}\subset\R^{1,0}&\R^{1,1}\subset\R^{1,1}&\R^{1,0}\subset\R^{2,0}&\R^{1,0}\subset\R^{1,1}&\R^{1,1}\subset\R^{2,2}\\

\hline
\hline
\multicolumn{6}{|c|}{\text{In case {\bf (2$^{\C}$)}, viewed as a $\underline J$-complex algebra, $\goth{s}=\langle L\rangle$ is given as in {\bf (2')}.}}\\
\hline
\end{array}$\hspace*{\fill}\\[1ex]
\hspace*{\fill}$\begin{array}{|c|ccc|}
\multicolumn{1}{c}{}&\multicolumn{1}{c}{\text{\bf (3)}}&\multicolumn{1}{c}{\text{\bf (3$^{\prime}$)}}&\multicolumn{1}{c}{\text{\bf (3$^{\C}$)}}\\
 \hline
\goth s=&\begin{nnarray}{c}\langle J_1,J_2,J_3\,|\\J_{[i]}J_{[i+1]}=J_{[i+2]},\\J_iJ_{i'}=-J_{i'}J_{i}\rangle\end{nnarray}& \begin{nnarray}{c}{\langle L_1,L_2,J\,|}\\J=-L_1L_2=L_2L_1,\\L_1=L_2J=-JL_2,\\L_2=JL_1=-L_1J\rangle\end{nnarray}&\begin{nnarray}{c}{\langle \underline J,L_1,L_2,J\,|}\underline J\in Z(\goth{s}),\\J=-L_1L_2=L_2L_1,\\L_1=L_2J=-JL_2,\\L_2=JL_1=-L_1J\rangle\end{nnarray}\\
\hline
\goth s\simeq&\HH&\mat_2(\R)&\mat_2(\C)\\
\hline
\begin{nnarray}{c}\text{through the}\\\text{morphism}\\\text{$\varphi$ given by:}\end{nnarray}&\begin{nnarray}{c}(J_1,J_2,J_3)\mapsto\\\text{the canonical}\\(\mathrm{i,j,k})\subset\HH\end{nnarray}
&\begin{nnarray}{c}(L_1,L_2,J)\mapsto\mbox{\footnotesize$\biggl($}
\mbox{\footnotesize$\left(\!\begin{narray}{cc}1&0\\0&-1\end{narray}\!\right)$},\\ 
\mbox{\footnotesize$\left(\!\begin{narray}{cc}0&1\\1&0\end{narray}\!\right)$},
\mbox{\footnotesize$\left(\!\begin{narray}{cc}0&-1\\1&0\end{narray}\!\right)$}\mbox{\footnotesize$\biggr)$}\end{nnarray}
&\begin{nnarray}{c}(\underline J,L_1,L_2,J)\mapsto\mbox{\footnotesize$\biggl($}{\mathrm i}I_2,\\
\mbox{\footnotesize$\left(\!\begin{narray}{cc}1&0\\0&-1\end{narray}\!\right)$},
\mbox{\footnotesize$\left(\!\begin{narray}{cc}0&1\\1&0\end{narray}\!\right)$},
\mbox{\footnotesize$\left(\!\begin{narray}{cc}0&-1\\1&0\end{narray}\!\right)$}\mbox{\footnotesize$\biggr)$}\end{nnarray}\\
\hline
\begin{nnarray}{c}\text{$\varphi$ conjugates}\\\text{adjunction}\\\text{in $\goth s$ to:}\end{nnarray}&\begin{nnarray}{c}\text{quaternionic}\\\text{conjugation}\\z\mapsto \overline{z}\end{nnarray}&
\begin{nnarray}{c}\mbox{\footnotesize$\left(\!\begin{narray}{cc}a&b\\c&d\end{narray}\!\right)\mapsto\left(\!\begin{narray}{cc}d&-b\\-c&a\end{narray}\!\right)$}\\\text{i.e. transpose of}\\\text{ the comatrix}\end{nnarray}
&\text{idem}\\
\hline
\varphi(\goth s^+)=&\R&\R.I_2&\C.I_2\\
\hline
\begin{nnarray}{c}\varphi(\goth s^+)\subset\varphi(\goth s)\\[-1ex]\text{\scriptsize[1]}\end{nnarray}&\R^{1,0}\subset\R^{4,0}&\R^{1,0}\subset\R^{2,2}&\R^{1,1}\subset\R^{4,4}\\
\hline
\hline
\multicolumn{4}{|c|}{\begin{nnarray}{c}\text{In case {\bf (3$^{\C}$)}, viewed as a $\underline J$-complex algebra,}\\{\goth s}=\langle L_1,L_2,L_3\,|\mathrm{i}L_{[i]}L_{[i+1]}=L_{[i+2]},L_{[i]}L_{[i']}=-L_{[i']}L_{[i]}\rangle.\end{nnarray}}\\
\hline
\end{array}$\hspace*{\fill}\\[.5ex]
$^\text{[1]}$given as (pseudo-)euclidean spaces for $\langle\cdot,\cdot\rangle$ introduced in Remark \ref{produit_naturel}; $\R^{a,b}$ means $(R^{a+b},\langle\cdot,\cdot\rangle)$ with sign$(\langle\cdot,\cdot\rangle)=(a,b)$.

\caption{\label{table_1}{\em Theorem \ref{structure_s} summed up in a table.} The first line gives $\goth s$ as a unital $\R$-algebra generated by $\langle\cdot,\cdot\rangle$-orthogonal complex and paracomplex structures. All letters $J$ denote complex structures, and $L$ paracomplex ones. All are $g$-skew adjoint, except the $g$-self adjoint underlined $\underline J$. Bracketed indices $[i]$ denote indices modulo $3$; $\goth s^+$ stands for the subspace of self adjoint elements of $\goth s$.
}
\end{table}

\begin{table}[ht!]
\hspace*{\fill}$\begin{array}{ccccccccc}
\multicolumn{1}{c}{}&\multicolumn{1}{c}{\text{\bf (1)}}&\multicolumn{1}{c}{\text{\bf (1$^{\C}$)}}&\multicolumn{1}{c}{\text{\bf (2)}}&\multicolumn{1}{c}{\text{\bf (2$^{\prime}$)}}&\multicolumn{1}{c}{\text{\bf (2$^{\C}$)}}&\multicolumn{1}{c}{\text{\bf (3)}}&\multicolumn{1}{c}{\text{\bf (3$^{\prime}$)}}&\multicolumn{1}{c}{\text{\bf (3$^{\C}$)}}\\
\hline
\begin{nnarray}{c}
\text{Condition}\\\text{on $\dim\MM$}\end{nnarray}&\begin{nnarray}{c}\text{any}\\\text{value}\end{nnarray}&\text{even}&\text{even}&\text{even}&\begin{nnarray}{c}\text{divisible}\\\text{by 4}\end{nnarray}&\begin{nnarray}{c}\text{divisible}\\\text{by 4}\end{nnarray}&\begin{nnarray}{c}\text{divisible}\\\text{by 4}\end{nnarray}&\begin{nnarray}{c}\text{divisible}\\\text{by 8}\end{nnarray}\\
\hline
\begin{nnarray}{c}
\text{Possible}\\\text{value of}\\\sign(g)\end{nnarray}&\begin{nnarray}{c}\text{any}\\\text{value}\end{nnarray}&\begin{nnarray}{c}(p,p)\\\text{with}\\p\in\N^\ast\end{nnarray}&\begin{nnarray}{c}(2p,2q)\\\text{with}\\p,q\in\N^\ast\end{nnarray}&\begin{nnarray}{c}(p,p)\\\text{with}\\p\in\N^\ast\end{nnarray}&\begin{nnarray}{c}(2p,2p)\\\text{with}\\p\in\N^\ast\end{nnarray}&\begin{nnarray}{c}(4p,4q)\\\text{with}\\p,q\in\N^\ast\end{nnarray}&\begin{nnarray}{c}(2p,2p)\\\text{with}\\p\in\N^\ast\end{nnarray}&\begin{nnarray}{c}(4p,4p)\\\text{with}\\p\in\N^\ast\end{nnarray}\\
\hline
\multicolumn{9}{c}{\text{\footnotesize(Of course, if $\sign(g)=(r,s)$ then $r+s=\dim\MM$.)}}\\
\end{array}$\hspace*{\fill}\\[1ex]
\hspace*{\fill}$\begin{array}{c|l|}
\multicolumn{1}{c}{}&\multicolumn{1}{c}{\text{In a well-chosen basis, with $p$ and $q$ the integers given above:}}\\
\hhline{~-}
\text{\bf (1)}&\mathrm{Mat}(g)=I_{p,q}\\
\hhline{~-}
\text{\bf (1$^{\C}$)}&\begin{nnarray}{l}\mathrm{Mat}(g)=I_{p,p},\ \mathrm{Mat}(\underline{J})=J_{p}\end{nnarray}\\
\hhline{~-}
\text{\bf (2)}&\begin{nnarray}{l}\mathrm{Mat}(g)=\diag(I_{p,q},I_{p,q}),\ \mathrm{Mat}(J)=J_{d/2}\end{nnarray}\\
\hhline{~-}
\text{\bf (2')}&\begin{nnarray}{l}\mathrm{Mat}(g)=I_{p,p},\ \mathrm{Mat}(L)=L_{p}\\
\text{or {\em e.g.\@}: }\mathrm{Mat}(g)=L_{p},\ \mathrm{Mat}(L)=I_{p,p}.\\\end{nnarray}\\
\hhline{~-}
\text{\bf (2$^{\C}$)}&\begin{nnarray}{l}\mathrm{Mat}(g)=L_{2p},\ \mathrm{Mat}(L)=I_{2p,2p},\ \mathrm{Mat}(\underline{J})=\diag(J_{p},-J_{p})\end{nnarray}\\
\hhline{~-}
\text{\bf (3)}&\begin{nnarray}{l}\mathrm{Mat}(g)=\diag(I_{p,q},I_{p,q},I_{p,q},I_{p,q}),\ \mathrm{Mat}(J_2)=J_{d/2},\\\mathrm{Mat}(J_1)=\diag(-J_{d/4},J_{d/4}),\ \mathrm{Mat}(J_3)=\mbox{\small$\left(\begin{narray}{cc}0&J_{d/4}\\J_{d/4}&0\end{narray}\right)$}\end{nnarray}\\
\hhline{~-}
\text{\bf (3')}&\begin{nnarray}{l}\mathrm{Mat}(g)=I_{2p,2p},\ \mathrm{Mat}(L_1)=L_{2p}\\\mathrm{Mat}({J})=\diag(-J_{p},J_{p})\ ^{\text{\bf [1]}},\ \mathrm{Mat}(L_2)=\mbox{\small$\left(\begin{narray}{cc}0&-J_{p}\\J_{p}&0\end{narray}\right)$}\end{nnarray}\\
\hhline{~-}
\text{\bf (3$^{\C}$)}&\begin{nnarray}{l}\mathrm{Mat}(g)=\diag(I_{2p,2p},-I_{2p,2p}),\ \mathrm{Mat}(\underline J)=\diag(J_{2p},J_{2p}),\\\mathrm{Mat}({J})=\diag(J_{p},J_{p},-J_{p},-J_{p}),\ \mathrm{Mat}(L_1)=L_{4p}\\\mathrm{Mat}(L_2)=\mbox{\small$\left(\begin{narray}{cc}0\\-\diag(J_{p},J_{p})\end{narray}\right.\hspace*{-0.2cm}\left.\begin{narray}{cc}\diag(J_{p},J_{p})\\0\end{narray}\right)$}\end{nnarray}\\
\hhline{~-}
\multicolumn{1}{l}{}&\multicolumn{1}{l}{^{\text{\bf [1]}}\,\text{or {\em e.g.\@}: }\mathrm{Mat}(g)=L_{2p},\ \mathrm{Mat}(L_1)=I_{2p,2p},\ \mathrm{Mat}({J})=\mbox{\small$\left(\begin{narray}{cc}0&J_{p}\\J_{p}&0\end{narray}\right)$}.}\\
\end{array}$\hspace*{\fill}\\[1.5ex]
\hspace*{\fill}$\begin{array}{|c|cccccccc|}
\multicolumn{1}{c}{}&\multicolumn{1}{c}{\text{\bf (1)}}&\text{\bf (1$^{\C}$)}&\text{\bf (2)}&\text{\bf (2')}&\text{\bf (2$^{\C}$)}&\text{\bf (3)}&\text{\bf (3')}&\multicolumn{1}{c}{\text{\bf (3$^{\C}$)}}\\
\hline
\begin{nnarray}{c}\text{ ``Complex Riemannian'' }\\\text{structures}\end{nnarray}&\varnothing&\{\pm\underline J\}&\varnothing&\varnothing&\{\pm\underline J\}&\varnothing&\varnothing&\{\pm\underline J\}\\
\hline
\text{K\"ahler structures}&\varnothing&\varnothing&\{\pm J\}&\varnothing&\{\pm J\}&\text{\scriptsize\bf [3]}&\text{\scriptsize\bf [4]}&\text{\scriptsize\bf [6]}\\
\hline
\text{para K\"ahler structures}&\varnothing&\varnothing&\varnothing&\{\pm L\}&\{\pm L\}&\varnothing&\text{\scriptsize\bf [5]}&\text{\scriptsize\bf [7]}\\
\hline
\hline
\multicolumn{9}{|l|}{%
^{\text{\bf [3]}}\text{the 2-sphere }\{\alpha J_1+\beta J_2+\gamma J_3;\alpha^2+\beta^2+\gamma^2=1\}=\{U;\langle U,U\rangle=1\}}\\
\multicolumn{9}{|l|}{%
^{\text{\bf [4]}}\text{the two-sheet hyperboloid }
\{\alpha L_1+\beta L_2+\gamma J;\alpha^2+\beta^2-\gamma^2=-1\}=\{U;\langle U,U\rangle=1\}
}\\
\multicolumn{9}{|l|}{%
^{\text{\bf [5]}}\text{the one-sheet hyperboloid}
\{\alpha L_1+\beta L_2+\gamma J;\alpha^2+\beta^2-\gamma^2=1\}=\{U;\langle U,U\rangle=-1\}
}\\
\multicolumn{9}{|l|}{%
^{\text{\bf [6]}}\text{the proper quadric with centre }\{U;\langle U,U\rangle_{\underline J}=1\}=
}\\
\multicolumn{9}{|l|}{%
\{\alpha L_1+\beta L_2+\gamma J;\alpha=\alpha'+\alpha''\underline J\ \text{\em etc. }\alpha^2+\beta^2-\gamma^2=-1\}}\\
\multicolumn{9}{|l|}{%
^{\text{\bf [7]}}\text{the proper quadric with centre }\{U;\langle U,U\rangle_{\underline J}=-1\}=
}\\
\multicolumn{9}{|l|}{%
\{\alpha L_1+\beta L_2+\gamma J;\alpha=\alpha'+\alpha''\underline J\ \text{\em etc. }\alpha^2+\beta^2-\gamma^2=1\}}\\
\hline
\end{array}$\hspace*{\fill}
\caption{\label{table_2}{\em Th.\@ \ref{structure_s} in a matricial form, with the sets of all (para)complex structures.}}
\end{table}

\begin{lem}\label{facteur_commun}
Take $U\in\goth{e}$ and $N\in\goth n$, $\mu$ the minimal polynomial of $U$ and $\mu'$ that of $U+N$. Then any irreducible factor of $\mu$ is also in $\mu'$, and {\em vice versa}.
\end{lem}

\noindent{\bf Proof.} $\mu(U+N)=\mu(U)+NV=NV$ with $V$ some polynomial in $U$ and $N$. As $\goth{n}$ is an ideal, $NV\in\goth n$ and by Proposition \ref{n_nilpotent}, $NV$ is nilpotent. So for some $k\in\N$, $(\mu^k)(U+N)=0$ {\em i.e.\@} $\mu'|\mu^k$. Symmetrically, $\exists l\in\N^\ast:\mu|\mu^{\prime l}$, so the result.\hfill{\rm $\Box$}\medskip

\noindent{\bf Proof of Theorem \ref{structure_s}.} We denote $T_m\MM$ by $E$ in this proof. We first state the announced classical results in associative algebra.

\begin{te}\label{wa} {\bf [Wedderburn -- Artin]} (see \cite{jacobson_structureofrings} \S3, p.\@ 40). Let $A$ be a finite dimensional semi-simple $\R$-algebra. Then $A$ is isomorphic to a direct sum of matrix algebras: $$A\simeq\underset{i=1}{\overset k\oplus}\mathrm{M}_{d_i}(\K_i)$$with for each $i$, $d_i\in\N^\ast$ and $\K_i\in\{\R,\C,\HH\}$.
\end{te}

\begin{te}\label{sn} {\bf [corollary of a theorem of Skolem -- Noether]} (see \cite{bourbaki_algebre8} \S10, no.\@ 1). Let $\theta$ be an automorphism of a finite dimensional semi-simple $\R$-algebra $A$. If $\theta$ acts trivially on the center of $A$, $\theta$ is interior.
\end{te}

As $g$ is orthogonally indecomposable, then if it is flat, $\dim\MM=1$ and $\goth s=\goth e=\R\Id$ is of type {\bf (1)}. We now suppose that $g$ is not flat. The list \ref{structure_s} follows from the orthogonal indecomposability of the action of $\goth{h}$, through the claim below. The elimination of only one possible form for $\goth{s}$ will also require, through Proposition \ref{pseudocommutation}, the fact that $\goth{h}$ is a holonomy algebra {\em i.e.\@} from the Bianchi identity satisfied by the curvature tensor.\\

\noindent{\em Claim \refstepcounter{claimcounter}\label{proj}\arabic{claimcounter}.}\label{polynome_minimal} If $p\in\goth{e}$ is self adjoint, its minimal polynomial $\mu_p(X)$ is of the form $Q^\alpha$ with $Q$ irreducible ---~so if $p$ is not invertible, it is nilpotent.\medskip

\noindent{\em Proof.} The minimal polynomial reads $\mu_p(X)=\prod_{i=1}^NQ_i^{\alpha_i}$ with irreducible $Q_i$'s. As $p$ is self adjoint, the direct sum $E=\oplus_{i=1}^N\ker Q^{\alpha_i}(p)$ is $g$-orthogonal. As $p\in\End(E)^\goth{h}$, each $\ker Q^{\alpha_i}(p)$ is $\goth{h}$-stable. Now $E$ is indecomposable, so $N=1$ and the claim.\\

By \ref{wm}, $\goth{e}=\goth{s}\oplus\goth{n}$ where $\goth{n}=\Rad(\goth{e})$ and  $\goth{s}$ is a semi-simple, self adjoint subalgebra of $\goth{e}$. As $\goth{n}$ is the intersection of the maximal ideals of $\goth{e}$ and as the adjunction $\sigma$ is an anti-morphism, $\goth{n}$ is also self adjoint. So \ref{wa} gives an isomorphism $\varphi:\goth{s}\overset{\simeq}{\longrightarrow}A$ with $A=\oplus_{i=1}^{k}I_i$ and $I_i=\mathrm{M}_{\delta_i}(\K_i)$, $\K_i\in\{\R,\C,\HH\}$. By a slight abuse, we also denote by $\sigma$ the conjugate action $\varphi\circ\sigma\circ\varphi^{-1}$ of $\sigma$ on $A$.\\

\noindent{\em Claim \stepcounter{claimcounter}\arabic{claimcounter}.} $k\leqslant2$. If $k=2$, then $\sigma(I_1)=I_2$. We then denote the $\delta_i$ by $\delta$ and the $\K_i$ by $\K$.\medskip

\noindent{\em Proof.} Let us denote by $1$ the unit matrix of any factor of $A$. As an (anti) automorphism of $A$, $\sigma$ acts on the factors $I_i$ of $A$, permuting them. Take $p=(1,0,\ldots,0)\in A$. As $p^2=p$, $\varphi^{-1}(p)$ is a (non zero) projection, so by Claim \ref{proj}, either $\varphi^{-1}(p)=1_\goth{e}$ and thus $k=1$, or $\sigma(p)\neq p$. In the latter case, $\sigma(I_1)\neq I_1$. Take $p'=p+\sigma(p)$. It is self adjoint by construction, and $p'^2=p^2+\sigma(p)^2=p+\sigma(p^2)=p'$ so it is a (non zero) projection. By Claim \ref{proj}, $\varphi^{-1}(p')=1_\goth{e}$ so $A=I_1\oplus\sigma(I_1)$ and then $k=2$.\\

\noindent{\em Claim \refstepcounter{claimcounter}\arabic{claimcounter}\label{pasH+H}.} If $k=2$, then $\delta=1$ and $\K=\R$ or $\K=\C$.\medskip

\noindent{\em Proof.} Suppose $k=2$. Take $p=(\diag(1,0,\ldots,0),0)\in A$ and $p'=p+\sigma(p)$. By the same reasoning as above, $\varphi^{-1}(p')$ is a non zero self ajoint projection so $p'=1_A$ by Claim \ref{proj}. As the $I_1$-component of $\sigma(p)$ is zero, then in fact $p=(1,0)$ and $\sigma(p)=(0,1)$; in particular $\diag(1,0,\ldots,0)=1_{I_1}$ {\em i.e.\@} $\delta=1$. Now Proposition \ref{pseudocommutation} implies $\K\neq\HH$. Indeed, suppose that $\K=\HH$, denote by $\mathrm i$ and $\mathrm j$ two of the three canonical roots of $-1$ in $\HH$, take $m\in\MM$ and $x,y\in T_m\MM$. Set $I=\varphi^{-1}(\mathrm i,0)$ and $J=\varphi^{-1}(\mathrm j,0)$ in $\goth{e}=\varphi^{-1}(\HH\oplus\HH)$. Notice that the $I_1$ component of $\sigma((\mathrm i,0))$ is zero, so $I^\ast J=0$, similarly $IJ^\ast=0$. By construction, $I+I^\ast$ is self adjoint, so:
\begin{align*}
R(x,y).(I+{I^\ast})(J+{J^\ast})&=R(x,y).(J+{J^\ast})(I+{I^\ast})\ \text{by Proposition \ref{pseudocommutation},}\\
&=R(x,y).(JI+{J^\ast}{I^\ast})\\
&=-R(x,y).(IJ+{I^\ast}{J^\ast})\qquad\text{as in $\HH$, $\mathrm j\mathrm i=-\mathrm i\mathrm j$,}\\
&=-R(x,y).(I+{I^\ast})(J+{J^\ast}).
\end{align*}
So $R(x,y).(I+{I^\ast})(J+{J^\ast})=0$. Now $(I+{I^\ast})(J+{J^\ast})=IJ+(IJ)^\ast=\varphi^{-1}((\mathrm i\mathrm j,0)+\sigma((\mathrm i\mathrm j,0)))$ is invertible, so for any $m\in\MM$ and any $x,y\in T_m\MM$, $R(x,y)=0$ {\em i.e.\@} $(\MM,g)$ is flat, in contradiction with $\goth{s}\simeq\HH\oplus\HH$.\\

Let us suppose $k=1$ and finish the proof. Let $\tau$ be the transposition $u\mapsto ^t\!\!u$ in $A\simeq{\mathrm M}_\delta(\K)$, and $\overline{\tau}$ its composition $u\mapsto ^t\!\!\overline u$ with the conjugation, in case $\K\in\{\C,\HH\}$. Then for $\K\in\{\R,\C\}$, respectively $\K\in\{\C,\HH\}$, $\tau$, respectively $\overline\tau$, is an anti-morphism (of $\R$-algebra) of $A$. So either $\tau\circ\sigma$ or $\overline\tau\circ\sigma$ is an automorphism of $A$ and, for $\K\in\{\R,\HH\}$, it acts trivially on the center $Z(A)$ as $Z(A)=\K.I_\delta$. If $\K=\C$, either $\sigma\circ\tau$ or $\sigma\circ\overline\tau$ acts trivially on the center $Z(A)=\C.I_\delta$. Then Theorem \ref{sn} gives a $v\in A$ such that $\sigma:u\mapsto v^t\widetilde uv^{-1}$ with $\widetilde u=u$ if $\K=\R$, $\widetilde u=\overline u$ if $\K=\HH$ and $\widetilde u=u$ or $\widetilde u=\overline u$ if $\K=\C$. As $\sigma^2=\Id_\goth{e}$, $v^t\widetilde v^{-1}\in Z(A)$ {\em i.e.\@} $^t\widetilde v=\lambda v$ with $\lambda\in\R$ if $\K\in\{\R,\HH\}$ and $\lambda\in\C$ if $\K=\C$. Applying $\widetilde\tau$ on both sides, we get that $\lambda=\pm1$ (in the case $\K=\C$ and $\widetilde u=\overline u$, we get only $|\lambda|=1$, but replacing $v$ by an adequate element of $\C.v$ achieves even $\lambda=1$).

If we replace $\varphi$ by $c_w\circ\varphi$ with $c_w:u\mapsto w^{-1}uw$, then $v$ is replaced by $wv^t\widetilde w$ {\em i.e.\@} $v$ undergoes a basis change like the matrix of a  bilinear or $\widetilde{\ }$-sesquilinear form. So using a suitable $c_w$, and recalling that $^t\widetilde v=\lambda v$ with $\lambda=\pm1$, we may suppose:\medskip

-- in case $\lambda=1$, that $v=\diag(I_{\delta'},-I_{\delta''})$ with $\delta'+\delta''=\delta$ if $\K=\R$ or ($\K\in\{\C,\HH\}$ and $\widetilde u=\overline u$), and that $v=I_{\delta}$ if ($\K=\C$ and $\widetilde u=u$),\medskip

-- in case $\lambda=-1$, that $\delta$ is even and $v=\mbox{\footnotesize$\left(\begin{array}{cc}0&-I_{\delta/2}\\I_{\delta/2}&0\end{array}\right)$}$ if ($\K\in\{\R,\C\}$ and $\widetilde u=u$), and that $v=I_{\delta}.{\mathrm i}$ if $\K=\HH$.\medskip

Now all cases where $v$ is diagonal imply $\delta=1$. Indeed, if $v=\diag(I_{\delta'},-I_{\delta''})$ or $v=I_\delta$, set $p=\diag(1,0,\ldots,0)$, and if $\K=\HH$ and $v=I_{\delta}.{\mathrm i}$, set $p=\diag({\mathrm j},0,\ldots,0)$. Then $p$ is self adjoint, non nilpotent, so $p=1_A$ or $p=1_A.{\mathrm j}$ by Claim \ref{proj} {\em i.e.\@} $\delta=1$. So if $\delta\geqslant2$, then $\K\in\{\R,\C\}$, $\lambda=-1$, $\widetilde u=u$, $\delta$ is even and $v=\mbox{\footnotesize$\left(\begin{array}{cc}0&-I_{\delta/2}\\I_{\delta/2}&0\end{array}\right)$}$. Setting $p'=\diag(1,0,\ldots,0)\in{\mathrm M}_{\delta/2}(\K)$ we get $p=\diag(p',p')$ a self adjoint non nilpotent element of $A$, so $p$ is invertible by Claim \ref{proj} {\em i.e.\@} $\delta'=1$ {\em i.e.\@} $\delta=2$. So the only allowed cases are those listed in \ref{structure_s}:\medskip

-- if $k=1$ and $\delta=1$, $\K=\R$ and $\sigma:u\mapsto^t\!\!u=u$, or $\K=\C$ and $\sigma:u\mapsto^t\!\!u=u$, or $\K=\C$ and $\sigma:u\mapsto^t\!\!\overline u=\overline u$, or $\K=\HH$ and $\sigma:u\mapsto^t\!\!\overline u=\overline u$,\medskip

-- if $k=1$ and $\delta=2$, ($\K=\R$ or $\K=\C$) and $\sigma:u\mapsto v^tuv^{-1}$ with  $v=\mbox{\footnotesize$\left(\begin{array}{cc}0&-1\\1&0\end{array}\right)$}$
 {\em i.e.\@} $\sigma$ is as described in Table \ref{table_1},\medskip

-- if $k=2$ and $\delta=1$ {\em i.e.\@} $A=I_1\oplus I_2$ with $I_1\simeq I_2\simeq\K$, ($\K=\R$ or $\K=\C$) and $\sigma$ permutes $I_1$ and $I_2$. Composing possibly $\varphi$ with a suitable automorphism of $A$, we get simply $\sigma:(a,b)\mapsto(b,a)$.\medskip

The remaining informations given in Tables \ref{table_1} and \ref{table_2} follow from standard calculations. We give only the following details.

In Table \ref{table_1}, the given generators are a (pseudo\nobreakdash-)orthogonal family of $(\goth s,\langle\cdot,\cdot\rangle)$, indeed $\frac 1d\tr(L^\ast L)=\frac1d\tr(-L^2)=\frac1d\tr(-\Id)=-1$ or, in case {\bf (2$^{\C}$)}, $\frac 1d\tr(L^\ast J)=\frac1d\tr(-\underline J)=0$ as $\underline J$ is a complex structure. 

For the three last columns of Table \ref{table_2}, we must check that the different (para)complex structures $U$ announced are indeed the only ones. Notice that if $U\in\goth s^\asym$, $U^2=\pm\Id\Leftrightarrow U^\ast U=\mp\Id\Rightarrow\langle U,U\rangle=\mp1$.

In cases {\bf (3)}, {\bf (3')}, and {\bf (3$^{\C}$)}, after Proposition \ref{formes_stables}, the (pseudo\nobreakdash-)K\"ahler manifold $(\MM,g,J)$ admits a non zero complex volume form so is Ricci flat. See also another brief proof in Theorem \ref{ricci}.\medskip

Finally, in \S\ref{realisation} are built the (non-empty) sets of germs of metrics inducing each case, and Proposition \ref{holonomie_generique} and Remark \ref{generalite} show Remark \ref{types_pour_h} above and hence the last assertion of the theorem.\hfill{\rm $\Box$}

\begin{rem}\label{bianchi_indispensable} In Claim \ref{pasH+H} above, the use of the Bianchi identity, through Proposition \ref{pseudocommutation}, is necessary. Consider the case $E=\R^{8p}\simeq\HH^p\oplus\HH^p$ and $\mathrm H':=\{u\in\mathrm{GL}_p(\HH)^2\,:\,u=(u_1,^t\overline{u_1})\}\subset\mathrm{GL}_{8p}(\R)$. Then $H'$ preserves the non degenerate real quadratic form $(x_1,x_2)\mbox{$\mapsto^t\!\overline{x_1}$}.x_2$ on $E$, and its action is orthogonally indecomposable. Now $\goth{gl}(E)^{{\goth h}'}=(\Id_{\HH^p}.\HH)^2\subset\mathrm{GL}_p(\HH)^2\subset\mathrm{GL}_{8p}(\R)$ and thus $\goth{gl}(E)^{{\goth h}'}\simeq\HH\oplus\HH$.
\end{rem}

The following corollary of Theorem \ref{structure_s} may be noticed.

\begin{cor}
A metric $g$ admits parallel self adjoint complex structures exactly in cases {\bf (1$^\C$)}, {\bf (2$^\C$)} and {\bf (3$^\C$)}, and they are: $\{\pm \underline J+N; N\in\goth{n}_0\ \text{and}\ N\underline J=-\underline JN\}$.
\end{cor}

\noindent{\bf Proof.} Suppose that some $\underline J_0\in\goth e^{\sym}$ satisfies $\underline J_0^2=-\Id$. Take the decomposition $\underline J_0=S+N$ with $S\in\goth s^\sym$ and $N\in\goth n^\sym$. By Lemma \ref{facteur_commun}, and as the minimal polynomial of $\underline J_0$ is $X^2+1$, irreducible, $S^2=-\Id$ so we are in case {\bf (1$^\C$)}, {\bf (2$^{\C}$)} or {\bf (3$^{\C}$)} and $S=\pm\underline J$. Now $-\Id=\underline J_0^2=(\underline J+N)^2=-\Id+JN+NJ+N^2$. By Proposition \ref{pseudocommutation}, $JN-NJ\in\goth n_0$, so $N(2J+N)=JN+NJ+N^2-(JN-NJ)=-(JN-NJ)\in\goth n_0$. By Lemma \ref{facteur_commun}, $((2J+N)^2+4\Id)^k=0$ for some $k$, so $2J+N$ is invertible, so $N\in\goth n_0$, and as then $N^2=0$, $N\in\{U\in\goth n_0;JU=-UJ\}$.\hfill{\rm $\Box$}\medskip

Finally, it may be useful to list the different possible parallel tensors.

\begin{prop}\label{formes_stables}
In each case of Theorem \ref{structure_s}, the metric admits the nondegenerate parallel multi- or sesquilinear forms of Table \ref{table_3} p.\@ \pageref{table_3}.
\end{prop}

\begin{table}[ht!]
\hspace*{\fill}$\begin{narray}{|ccc|}
\hline
\text{parallel tensor/exists in cases}&\text{parametrised by}&\text{given as}\\\hline
\begin{nnarray}{c}
\text{Pseudo-Riemannian}\\\text{metric/all cases}
\end{nnarray}&U\in\goth e^{\sym}\smallsetminus\goth n^{\sym}&g(\,\cdot\,,U\cdot\,){}\\[2.5ex]
\begin{nnarray}{c}\text{Symplectic form/all}\\\text{except {\bf (1)} and {\bf (1$^\C$)}}\end{nnarray}&\begin{nnarray}{c}U=V+N,\\V\in(\goth s^{\asym})^\ast, N\in\goth n^{\asym}\end{nnarray}&g(\,\cdot\,,U\cdot\,){}\\[2.5ex]
\begin{nnarray}{c}
\text{``Complex Riemannian''}\\\text{metric/{\bf (1$^\C$)}, {\bf (2$^\C$)}, {\bf (3$^\C$)}}
\end{nnarray}&\begin{nnarray}{c}U\in\goth e^{\sym}\smallsetminus\goth n^{\sym}\\\text{such that }U\underline J=\underline JU\end{nnarray}&\begin{nnarray}{c}\underline g_U=\\g(\,\cdot\,,U\cdot\,)+{\rm i}g(\,\cdot\,,\underline JU\cdot\,)\end{nnarray}{}\\[2.5ex]
\begin{nnarray}{c}
\text{Hermitian (pseudo-)K\"ahler}\\\text{metric w.\@ r.\@ to some $J\in\goth s^{\asym}$}\\\text{{\bf (2)}, {\bf (2$^{\C}$)}, {\bf (3)}, {\bf (3')}, {\bf (3$^{\C}$)}}\end{nnarray}&\begin{nnarray}{c}U\in\goth e^{\sym}\smallsetminus\goth n^{\sym}\\\text{such that }UJ=JU\end{nnarray}&\begin{nnarray}{c}h_U=\\g(\cdot,U\cdot)+{\rm i}g(\cdot,JU\cdot)\end{nnarray}{}\\[3ex]
\begin{nnarray}{c}
\text{$\underline J$-complex symplectic form}\\\text{{\bf (2$^\C$)}, {\bf (3$^\C$)}}\end{nnarray}&\begin{nnarray}{c}U=V+N,\\V\in(\goth s^{\asym})^\ast, N\in\goth n^{\asym},\\\text{such that }N\underline J=\underline JN\end{nnarray}&\begin{nnarray}{c}\underline\omega_U=\\g(\cdot,U\cdot)+{\rm i}g(\cdot,\underline JU\cdot)\end{nnarray}{}\\[3.5ex]
\begin{nnarray}{c}\text{$J$-complex symplectic form}\\\text{{\bf (3)}, {\bf (3')}, {\bf (3$^{\C}$)}}\end{nnarray}&\begin{nnarray}{c}U=V+N,\\V\in(\goth s^{\asym})^\ast, N\in\goth n^{\asym},\\\text{such that }UJ=-JU\end{nnarray}&\begin{nnarray}{c}\omega_U=\\g(\cdot,U\cdot)+{\rm i}g(\cdot,JU\cdot)\end{nnarray}{}\\[2.5ex]
\begin{nnarray}{c}\text{Non null $\underline J$-complex volume}\\\text{form/{\bf (1$^\C$)}, {\bf (2$^\C$)}, {\bf (3$^\C$)}}\end{nnarray}&\multicolumn{2}{c|}{\text{associated with $\underline g_U$ above}}{}\\[2.5ex]
\begin{nnarray}{c}\text{Non null $J$-complex volume}\\\text{form/{\bf (3)}, {\bf (3')}, {\bf (3$^{\C}$)}}\end{nnarray}&\multicolumn{2}{c|}{\text{equal to }\omega_U^{\wedge(d/4)}\ \text{with $\omega_U$ as above}}\\\hline
\end{narray}$\hspace*{\fill}
\caption{\label{table_3}{\em The real and complex parallel tensors existing in the different cases.} In cases {\bf (3)}, {\bf (3')} and {\bf (3$^\C$)}, $(\goth s^{\asym})^\ast$ is the complement of the isotropic cone in $\goth s^{\asym}$. The real part of $h_U$ is a (pseudo\nobreakdash-)\-Riemannian metric, its imaginary part is a 2-form of type (1,1).}
\end{table}

\noindent{\bf Proof.} Some lines of Table \ref{table_3} require a brief checking.

{\bf (1)} Any $U\in\goth e^{\sym}\smallsetminus\goth n^{\sym}$ is nondegenerate. Indeed, any $U\in\goth s^{\sym}\smallsetminus\{0\}$ is (see Table \ref{table_1}), so its minimal polynomial $\mu$ is not divisible by $X$; by Lemma \ref{facteur_commun}, neither is the minimal polynomial of $U+N$ for any $N\in\goth n$.

{\bf (2)} If some nondegenerate alternate form is parallel for a torsion-free connection, it is closed, thus symplectic. Then proceed as in {\bf (1)} above.

{\bf (3)} If $J$ is a parallel complex structure (self- or skew-adjoint), nondegenerate complex bilinear forms are the $g(\cdot,U\cdot)-{\mathrm i}g(\cdot,V\cdot)$ such that (check it) $\ker U\cap\ker V=\{0\}$, $V=UJ$, $U^\ast=U$ and $V^\ast=V$. By Proposition \ref{n_nilpotent}, the first condition implies that $U\not\in\goth n$ or $V\not\in\goth n$, so by Lemma \ref{facteur_commun} and the reasoning of {\bf (1)}, that $U$ or $V$ is nondegenerate, hence both. Now if $J^\ast=-J$, the relations give that $UJ=-JU$ and $VJ=-JV$. As $U^\ast=U$, by Proposition \ref{pseudocommutation}, everywhere, $R(\,\cdot\,,\,\cdot\,)(UJ-JU)=0$. As $UJ-JU=2UJ$ is nondegenerate, $\MM$ would be flat. So $J^\ast=J$, we denote it by $\underline J$. This time $U\underline J=\underline JU$. After Table \ref{table_1} and Lemma \ref{facteur_commun}, the existence of such a $\underline J$ leads to the announced form of $\goth s^{\sym}$. Conclude by the same reasoning as in {\bf (1)}; {\bf (4)}-{\bf (6)} are entirely similar.

{\bf (7)} If some parallel $\underline J$ exists, so some complex Riemannian metric $\underline g_{\underline J}$ as on line 3 of Table \ref{table_3}, take $(e_i)_{i=1}^{d/2}$ some $\underline g_{\underline J}$-orthonormal complex frame field, and $\nu=e_1^\ast\wedge\ldots\wedge e_{d/2}^\ast$. As $\underline g_{\underline J}$ is parallel, so is $\nu$.\hfill{\rm $\Box$}

\section{The space of germs of metrics realising each form of~$\goth s$}\label{realisation}

\begin{reminder}\label{rappels_constructions} Metrics with $\goth s$ of type {\bf (1$^\C$)} are the real parts of complex Riemannian metrics {\em i.e.\@} of holomorphic, non degenerate $\C$-bilinear forms on complex manifolds $(\MM,\underline J)$. It is well known and easy to check.

As it is also well known, germs of (pseudo\nobreakdash-)\-K\"ahler metrics (type {\bf (2)}) are parametrised by a K\"ahler potential $u$, which is a real function:\medskip\\
\phantom{{\bf (a)}}\hspace*{\fill}$\displaystyle g\left(\frac{\partial}{\partial z_i},\frac{\partial}{\partial \zbar_j}\right)=\frac{\partial^2u}{\partial z_i\partial\overline z_j}$.\hspace*{\fill}{\bf (a)}\medskip\\
Similarly, germs of para K\"ahler metrics (type {\bf (2')}) are parametrised by a para K\"ahler potential (see {\em e.g.\@} \S2 of \cite{BB-ikemakhen1997}). The supplementary distributions $\ker(L\pm\Id)$ 
are integrable. Take $((x_i)_{i=1}^{d/2},(y_i)_{i=1}^{d/2})$ coordinates adapted to the corresponding pair of integral ($g$-isotropic) foliations. Then the metrics of type {\bf (2')} depend on a real potential $u$ through:\medskip\\
\phantom{{\bf (b)}}\hspace*{\fill}$\displaystyle g\left(\frac{\partial}{\partial x_i},\frac{\partial}{\partial y_j}\right)=\frac{\partial^2u}{\partial x_i\partial y_j}$.\hspace*{\fill}{\bf (b)}\medskip\\
A metric of type {\bf (2$^\C$)} is given by the complexification of {\bf (a)} or {\bf (b)}, indifferently: take $u$ complex and replace the real and imaginary parts of the $z_i$, in case {\bf (a)}, or $(x_i)_i$ and  $(y_i)_i$, case {\bf (b)}, by complex variables.
\end{reminder}

\begin{rem}Be careful however that a manifold $(\MM,g)$ of type {\bf (2)} or {\bf (2$^\C$)} has to be complex, hence in particular real analytic, whereas one of type {\bf (2')} may be only smooth.
\end{rem}

\begin{rem}We recall also that the ``complex Riemannian'' metrics defined in Table \ref{table_3} in cases {\bf (1$^\C$)}, {\bf (2$^\C$)} and {\bf (3$^\C$)} are holomorphic with respect to the self adjoint complex structure $\underline J$. Check that, if $z_j=x_j+{\rm i}y_j$ are complex coordinates, $\frac{\partial}{\partial y_j}\underline{g}_{k,l}={\rm i}\frac{\partial}{\partial x_j}\underline{g}_{k,l}$ for all $k,l$.
\end{rem}

\begin{prop}\label{holonomie_generique}A generic metric of type {\bf (1)}, {\bf (2)}, {\bf (2')}, {\bf (1$^\C$)} or {\bf (2$^\C$)} has the holonomy algebra given in Remark \ref{types_pour_h}. More precisely, if the 2-jet at the origin of some metric of the wished type satisfies some dense open condition among such 2-jets, then its holonomy algebra is as in Remark \ref{types_pour_h}. In particular, those holonomy groups are obtained on a dense open subset, for the $C^2$ topology, of the corresponding metrics.
\end{prop}

\noindent{\bf Proof.} It is standard, but we did not find any really explicit reference in the literature, and we need such a reference, as we will generalise it in our work on $\goth n$.  Besides it is short, and makes this paper self-contained. So we recall it. At the origin, take normal coordinate vectors $(X_i)_{i=1}^d$, moreover such that $X_{i+1}=JX_i$ or $X_{i+1}=LX_i$ for $i$ odd, in case {\bf (2)} or {\bf (2')}. So for any coordinate vectors $U,V$, $D_UV=0$ at $0$. For any coordinate vectors $A,B,U,V$ at the origin, $g(R(A,B)U,V)$ is equal to (check it):
\begin{align*}
&\frac12\bigl(A.U.(g(B,V))-B.U.(g(A,V))-A.V.(g(B,U))+B.V.(g(A,U))\bigr).
\end{align*}
In case {\bf (1)}, $g(R(X_i,X_j)_{|0}\,\cdot\,,\,\cdot\,)$ is the alternate part of the bilinear form:
$$\beta_{i,j}:(U,V)\mapsto X_i.U.g(X_j,V)-X_j.U.g(X_i,V).$$
The $\beta_{i,j}$ depend on the second derivatives of the coefficients of $g$ at $0$, which are free in normal coordinates. So, on a dense open subset of the 2-jets of metrics, their alternate parts are linearly independent and span a $\frac{d(d-1)}2$-dimensional space in $\dim\goth{o}_d(\R)$ {\em i.e.\@} $\goth{o}_d(\R)$ itself.

In case {\bf (2)} we set, for $j$ odd, $Z_{\frac{j+1}2}=X_j-{\rm i}X_{j+1}$ and $\Zbar_{\frac{j+1}2}=X_j+{\rm i}X_{j+1}$ in $T^\C\MM$. The $R(Z_i,Z_j)$ and $R(\Zbar_i,\Zbar_j)$ vanish, and the $R(Z_i,\Zbar_j)$ vanish when evaluated on $\Lambda^2T^{1,0}\MM$ or $\Lambda^2T^{0,1}\MM$. So $R$ is determined at $0$ by the $\beta_{i,j}$: $(Z_k,\Zbar_l)\mapsto\ g(R(Z_i,\Zbar_j),Z_k,\Zbar_l)$. As:
$$g(R(Z_i,\Zbar_j),Z_k,\Zbar_l)={\textstyle\frac12}(-\Zbar_j.Z_k.(g(Z_i,\Zbar_l))-Z_i.\Zbar_l.(g(\Zbar_j,Z_k))),$$
$R_{|0}$ is given by the fourth derivatives of the K\"ahler potential $u$.  Those are free in normal coordinates, so on a dense open subset of the 2-jets of metrics, the $(\beta_{i,j})_{i,j=1}^{d/2}$ are linearly independent hence span a $\left(\frac{d}2\right)^2$-dimensional space in $\goth{u}_{d/2}$, hence $\goth{u}_{d/2}$.

For {\bf (2')}, replace $(Z_i,\Zbar_i)_{i=1}^{d/2}$ by $(X_{i},Y_{i})_{i=1}^{d/2}$, and $\goth{u}_{d/2}$ by $\goth{gl}_{d/2}(\R)$.

For types {\bf (1$^\C$)} and  {\bf (2$^\C$)}, $R$ is $\underline J$-complex; repeat the proofs in complex coordinates.\hfill{\rm $\Box$}\medskip

Now we describe the space of germs of metrics of type {\bf (3)}, {\bf (3')} and {\bf (3$^\C$)}. It is classical for type {\bf (3)} (hyperk\"ahler), the other cases are an adaptation of the argument.

\begin{notation}
Take $\varepsilon\in\{-1,1\}$ and $\delta\in\N^\ast$. We denote by $\GG_\varepsilon$ the set of germs at 0 of triples $(g,J,U)$ with $g$ a (pseudo\nobreakdash-)\-Riemannian metric on $\R^d=\R^{4\delta}$ and $J$ and $U$ two $g$-skew adjoint parallel endomorphisms fields such that $\varepsilon U^2=-J^2=\Id$, and $JU=-UJ$. We define $\GG_\C$ similarly, with $g$ a complex Riemannian metric on $\C^{4 \delta}$ and similar $J$ and $U$ (with {\em e.g.\@} $\varepsilon=-1$, but this makes no difference on $\C$).
\end{notation}

Using Cartan-K\"ahler theory (see \cite{bcggg,ivey-landsberg}), we parametrise $\GG_\varepsilon$ and $\GG_\C$ in the real analytic category. We proceed as R.\@ Bryant did in \cite{bryant1996} \S2.5 pp.\@ 122--126 for hyperk\"ahler metrics {\em i.e.\@} for $\varepsilon=-1$, detailing the calculations to show that the case $\varepsilon=1$ works alike, and to allow another generalisation of this construction in our work on $\goth n$. The complex case $\GG_\C$ follows. This provides in particular an explicit writing of R.\@ Bryant's line of proof given in \cite{bryant1996}; we did not find this in the literature.

\begin{rem-notation}
Let $\omega_0$ be some complex symplectic form on some open set $\OO$ of $\C^{2\delta}$. Then any 2-form $\omega$ of type (1,1), real, may be written as $\omega=\Im\left(\omega_0(\,\cdot\,,U_\omega\,\cdot\,)\right)$, with $U_\omega$ an $\omega_0$-self adjoint complex antimorphism field. The correspondence is bijective between such forms $\omega$ and such $U_\omega$, so we use this notation $U_\omega$ in the following.
\end{rem-notation}

\begin{rem}\label{correspondence1}
 The set $\GG_\varepsilon$ is in bijection with the set $\GG'_\varepsilon$ of germs of couples $(\omega_0,\omega)$, with $\omega_0$ a complex symplectic form on $\C^{2\delta}$ and $\omega$ a closed 2-form of type (1,1), real, such that $U_\omega^2=\varepsilon\Id$, through the following.\medskip

-- Let $(g,J,U)$ be given. Then on $\C^{2\delta}:=(\R^{4\delta},J)$ set:
\begin{gather*}
\omega_0:=g(\,\cdot\,,U\,\cdot\,)+{\rm i}g(\,\cdot\,,JU\,\cdot\,)\ \text{ and }
\omega:=\varepsilon g(\,\cdot\,,J\,\cdot\,)=\Im\left(\omega_0(\,\cdot\,,U\,\cdot\,)\right).
\end{gather*}
As $DJ=DU=0$, immediately $\dd \omega_0=\dd\omega=0$.\medskip

-- Let $(\omega_0,\omega)$ be given. Then on $(\R^{4\delta},J):=(\C^{2\delta},{\rm i})$ set:
\begin{gather*}
g:=-\varepsilon\omega(\,\cdot\,,{\rm i}\,\cdot\,)\quad \text{and}\quad U:=U_\omega.
\end{gather*}
As $\dd \omega_0=\dd\omega=0$, $DJ=DU=0$. This is standard, see {\em e.g.\@} \cite{moroianu} \S 11.2.\medskip
\end{rem}

In this new point of view, up to a biholomorphism of $\C^{2\delta}$, $\omega_0$ may be considered, by the Darboux theorem, as the canonical symplectic form:
\begin{gather*}
\omega_0=\sum_{j=1}^{\delta}\dd z_i\wedge\dd z_{\delta +i}=\frac12{}^t\!\dd z\wedge\Omega_0\wedge\dd z\ \text{with }\Omega_0=\text{\small$\left(\begin{array}{cc}0&I_\delta\\-I_ \delta&0\end{array}\right)$},
\end{gather*}$\dd z$ denoting the column $(\dd z_i)_{i=1}^{2\delta}$. {\em From now on, we consider that $\omega_0$ is this canonical form.} Then the elements of $\GG_\varepsilon$, seen up to diffeomorphism of $\R^{4\delta}$, are in bijection with those of $\GG'_\varepsilon$, seen up to symplectomorphism of $(\C^{2\delta},\omega_0)$. Now we use Cartan-K\"ahler theory to describe $\GG'_\varepsilon$.

\begin{notation} Set $V:=\Mat(U)$, $U$ is an antimorphism so $U(z)=V.\zbar$. As $\omega_0(U\,\cdot\,,\,\cdot\,)=-\overline{\omega_0(\,\cdot\,,U\,\cdot\,)}$, we get $\Omega_0V=-\,^t\overline V\Omega_0$. A 2-form $\omega$ is in $\GG'_\varepsilon$ if and only if it is closed and:
\begin{gather*}
\omega=\Im(\omega_0(\,\cdot\,,U\,\cdot\,))={\textstyle \frac1{2{\rm i}}}{}^t\!\dd z\wedge\Omega_0V\wedge\dd \zbar\ \text{with }V\overline{V}=\varepsilon\Id
\end{gather*}
{\em i.e.\@}, setting $H:=-\Omega_0V$, if and only if:
\begin{gather*}
\omega={\textstyle\frac{\rm i}2}{}^t\!\dd z\wedge H\wedge\dd \zbar\ \text{with }^tH=\Hbar\ \text{and }\Hbar\Omega_0H=-\varepsilon\Omega_0.
\end{gather*}
Let $\HHH_\varepsilon\subset{\rm M}_{2\delta}(\C)$ be the space of such matrices $H$. The (1,1)-forms $\omega$ such that $U_\omega^2=\varepsilon\Id$ are exactly given by the functions $H:\C^{2\delta}\rightarrow\HHH_\varepsilon$, through: $\omega_H:=\frac{\rm i}2{}^t\!\dd z\wedge H(z)\wedge\dd \zbar$. Denoting by $(z,H)$ the points in $\C^{2\delta}\times\HHH_\varepsilon$, such an $\omega_H$ is {\em closed} if and only if the 3-form $\lambda:={}^t\!\dd z\wedge\dd H\wedge\dd\zbar$ vanishes along the graph $\SSS$ of $H$. So we are looking for the integral manifolds $\SSS$ of the exterior differential system ${\bf I}=(\lambda)$ on $\C^{2\delta}\times \HHH_\varepsilon$, with the independence condition that $\dd z_1\wedge\ldots\wedge\dd z_{2\delta}$ never vanishes ({\em i.e.\@} $\SSS$ is the graph of some $H:\C^{2\delta}\rightarrow\HHH_\varepsilon$). Then the Cartan-K\"ahler theorem parametrises $\GG'_\varepsilon$, hence $\GG_\varepsilon$, providing:
\end{notation}

\begin{prop}\label{parametrageG}
The elements of $\GG_\varepsilon$, considered up to diffeomorphism, are parametrised by $\frac d2=2\delta$ real analytic functions of $2\delta+1$ real variables. Those of $\GG_\C$, up to biholomorphism, are parametrised by $\frac d4=2\delta$ holomorphic functions of $2\delta+1$ complex variables.\end{prop}

\begin{rem}\label{parametrage3} The generality of the elements of $\GG_\varepsilon$ and
 $\GG_\C$ ensures that their corresponding  algebra $\goth s$ is indeed, generically, in cases {\bf (3)}, {\bf (3')} or {\bf (3$^\C$)} (and {\em e.g.\@} not the full $\End(T\MM)$). In fact, their holonomy group itself is generically that of Remark \ref{types_pour_h}, see Remark \ref{generalite}.
\end{rem}

\noindent{\bf Proof.} The writing of ${\bf I}$ in $\C^{2\delta}\times \HHH_\varepsilon$ does not depend on $z$, so we have only to perform Cartan's test on some arbitrary fibre $\{z_0\}\times H$, say with $z_0=0$. Moreover, over that point $z_0$, the symplectic group $\Sp(2\delta,\C)$ acts transitively on $\bigl\{\frac{\rm i}2{}^t\!\dd z\wedge H\wedge\dd \zbar ; H\in\HHH_\varepsilon\bigr\}$, preserving ${\bf I}$, so we have only to perform Cartan's test at some specific element $H_0\in\HHH_\varepsilon$, say:\medskip

-- if $\varepsilon=-1$, $H_0=I_{p,q,p,q}=\operatorname{diag}(I_p,-I_q,I_p,-I_q)$ with $p+q=n$,\medskip

-- if $\varepsilon=1$, $H_0={\rm i}I_{n,n}$.\medskip

\noindent{\em Remark.} As it appears in \cite{bryant1996}, the connected component $\HHH^{p,q}_{-1}$ of $I_{p,q,p,q}$ in $\HHH_{-1}=\sqcup_{p+q=n}\HHH^{p,q}_{-1}$ is canonically isomorphic to $\operatorname{Sp}(n,\C)/\operatorname{Sp}(p,q)$. So choosing some function $H:\C^{2\delta}\rightarrow\HHH^{p,q}_{-1}$ amounts to choosing a reduction to $\operatorname{Sp}(p,q)$, which is a real form of $\operatorname{Sp}(n,\C)$, of the principal bundle $\operatorname{Sp}(n,\C)\times\C^{2\delta}$. Similarly here, $\HHH_{1}\simeq\operatorname{Sp}(n,\C)/\operatorname{Sp}(n,\R)$ so choosing some $H:\C^{2\delta}\rightarrow\HHH_{1}$ is choosing a reduction of it to $\operatorname{Sp}(n,\R)$, which is another real form of $\operatorname{Sp}(n,\C)$.\medskip

Let us set $\partial z_j=\partial x_j+{\rm i}\partial y_j$. If a subspace $E$ of $T_{m_0}\MM$ is horizontal {\em i.e.\@} tangent to the factor $\C^{2\delta}$, $\lambda_{|E}=0$ so $E$ is an integral element of ${\bf I}$. Let us define $(E_k)_{k=0}^{4\delta}$ by:
\begin{gather*}
E_k=\Span\left((e_j)_{j=1}^k\right)\ \text{  with, for $1\leqslant j\leqslant \delta$: }\ e_j=\left(\partial x_j,0 \right)\text{ and}\\e_{\delta +j}=\left(\partial x_{\delta +j}+\frac{j-1}\delta\partial y_{\delta +j},0 \right),\ 
\text{and for $1\leqslant j\leqslant 2\delta$: }\ e_{2\delta+j}=\left(\partial y_{j},0\right).
\end{gather*}
Each $E_k$ is horizontal so $(E_k)_{k=0}^{4\delta}$ is an integral flag of ${\bf I}$ at $m_0$. We classically set $H(E_k):=\bigl\{v ; \Span(v,E_k)$ is an integral element of ${\bf I}\bigr\}$, and $s_k:=\codim_{H(E_{k-1})}H(E_{k})$ the $k$th. character of ${\bf I}$ (indeed this flag is ordinary, as we will see). We will check:\medskip

{\bf (1)} for all $k$, $s_k=k-1$, and $s_k=0$ for $k>2\delta+1$,\medskip

{\bf (2)} $\dim V_{4\delta}({\bf I})\geqslant 2C_{2\delta+2}^3$, with $V_{4\delta}({\bf I})$ the variety of integral elements of ${\bf I}$ in the grassmannian $G_{4\delta}(T(\C^{2\delta}\times\HHH_\varepsilon))$.\medskip

After Cartan's criterion, $\dim V_{4\delta}({\bf I})\leqslant\sum_{k=1}^{4\delta}ks_k$, and if equality holds then $E_{4\delta}$ is ordinary. So here:
$$\dim V_{4\delta}({\bf I})\leqslant\sum_{k=1}^{4\delta}ks_k=\sum_{k=1}^{2\delta+1}k(k-1)=\frac83\delta^3+4\delta^2+\frac43 \delta=2C_{2\delta+2}^3.$$
As $2C_{2\delta+2}^3\leqslant\dim V_{4\delta}({\bf I})$, equality holds, hence $E_{4\delta}$ is ordinary and after the Cartan-K\"ahler theorem, ${\bf I}$ admits an integral manifold $\SSS$ through $(0,H_0)$ with $T\SSS=E_{4\delta}$, and the space of germs of integral manifolds passing by $z_0$ depends on $s_{k}$ functions of $k$ variables, with $s_k$ the last non vanishing character, so here $2\delta$ functions of $2\delta+1$ variables. This parametrisation of the set $\GG_\varepsilon$ is done up to the choice of complex Darboux coordinates for $\omega_0$, and such coordinates depend on one symplectic generating function of $2\delta$ variables. As $2\delta<2\delta+1$, this does not interfer and $\GG_\varepsilon$ itself is parametrised by $2\delta$ functions of $2\delta+1$ variables, the proposition. We are left with showing {\bf (1)} and {\bf (2)}.\medskip

We introduce $W_\varepsilon:=T_{H_0}\HHH_{\varepsilon}$, then:
\begin{gather*}W_{1}=\left\{\text{\small$\left(\begin{array}{cc}a&b\\\overline{b}&\overline{a}\end{array}\right)$};a,b\in{\rm M}_ \delta(\C),\,^ta=\abar,\ ^tb=b\right\}\\
\text{and: }\  W_{-1}=\left\{\text{\small$\left(\begin{array}{cc}a&I_{p,q}b\\\overline{b}I_{p,q}&-I_{p,q}\overline{a}I_{p,q}\end{array}\right)$};a,b\in{\rm M}_n(\C),\,^ta=\abar,\ ^tb=b\right\}.
\end{gather*}
Then $T_{m_0}\MM_\varepsilon=T_0\C^{2\delta}\oplus W_\varepsilon\simeq\C^{2\delta}\oplus W_\varepsilon$ and the subset of the grassmannian $G_{4\delta}(T_{m_0}\MM)$ on which the independence condition holds is canonically identified with $(\C^{2\delta})^\ast\otimes W_\varepsilon$. 

{\bf (1)} follows from the fact that for $k>2\delta$, $H(E_{k})=\C^{2\delta}\oplus\{0\}$, and for $1\leqslant k\leqslant \delta$:\medskip

-- $H(E_k)=\C^{2\delta}\oplus\{\Im a_{i,j}=0$ for $1\leqslant i<j\leqslant k\}\subset\C^{2\delta}\oplus W_\varepsilon$, so $\codim_{H(E_{k-1})}H(E_k)=k-1$,
\medskip

-- $H(E_{n+k})=\C^{2\delta}\oplus\bigl\{\Re a_{i,j}=\Re b_{i,j}=0$ for $1\leqslant i<j\leqslant k$ and $\Im b_{k,j}+\frac{k-1}\delta\Re b_{k,j}=0$ for $k\leqslant j\leqslant \delta\bigr\}\subset\C^{2\delta}\oplus W_\varepsilon$, so $\codim_{H(E_{\delta+k-1})}H(E_{\delta +k})= \delta +k-1$.\medskip

To check {\bf (2)}, we introduce some notation. We denote the basis vectors $(\partial x_i)_{i=1}^{2\delta}$ of $\C^{2\delta}$ by $((u_i)_{i=1}^k,(u'_i)_{i=1}^k)$ (the $u_i$ and $u'_i$ are $\omega_0$-dual), then $(\partial y_i)_{i=1}^{2\delta}=((Ju_i)_{i=1}^k,(Ju'_i)_{i=1}^k)$. We denote by $H^{(1)}$ a generic element of $(\C^{2\delta})^\ast\otimes W_\varepsilon$. If a function $H:\C^{2\delta}\rightarrow\HHH_\varepsilon$ with $H(0)=H_0$ is such that $\dd H_{|0}=H^{(1)}$, then, at $0$, $\dd\omega_H$ is determined by $\dd\omega_H=\lambda_{|m_0}(H^{(1)}\,\cdot\,,H^{(1)}\,\cdot\,,H^{(1)}\,\cdot\,)$, that we denote by $\lambda_{H^{(1)}}$. In concrete terms, for the calculations below, $\lambda_{H^{(1)}}(u,v,w)$ is equal to:
$$\omega_0\left(u,H^{(1)}(v).\wbar\right)+\omega_0\left(v,H^{(1)}(w).\ubar\right)+\omega_0\left(w,H^{(1)}(u).\vbar\right).$$
At $(0,H_0)$, $V_{4\delta}({\bf I})$ is the set of the 1-jets of {\em closed} 2-forms $\omega_H$ as wanted. An $H^{(1)}$ is in $V_{4\delta}({\bf I})$ if and only if $\lambda_{H^{(1)}}=0$, which may be written as the two following conditions:
\begin{center}
\begin{tabular}{l}
{\bf (a)} for any three $\{u,v,w\}\subset\{u_i,Ju_i\}_{i=1}^k$, $\lambda_{H^{(1)}}(u^{(\prime)},v^{(\prime)},w^{(\prime)})=0$,\\[.5ex]
{\bf (b)} for any two $\{u,v\}\subset\{u_i,Ju_i\}_{i=1}^k$,
\\\hspace*{13em} $\lambda_{H^{(1)}}(u,u^{\prime},v^{(\prime)})=\lambda_{H^{(1)}}(v,v^{\prime},u^{(\prime)})=0$.
\end{tabular}
\end{center}
The parenthesised primes enable to denote several equations at once, so {\bf (a)} consists of $8C_{2\delta}^3$ equations and {\bf (b)} of $4C_{2\delta}^2$. Now the equations of {\bf (a)} are redundant. Indeed the reader may check the following. Take any $H^{(1)}$ and any $\{i,j,k\}\subset\llbracket1,\delta\rrbracket$ and $\{\alpha,\beta,\gamma\}\subset\{0,1\}$ such that $\sharp\{J^\alpha u_i,J^\beta u_j,J^\gamma u_k\}=3$ (so, $C_{2\delta}^3$ possibilities). Set $(u,v,w):=(J^\alpha u_i,J^\beta u_j,J^\gamma u_k)$ and, in case $\varepsilon=1$, $\eta_1:=(-1)^{\gamma-\beta}$, $\eta_2:=(-1)^{\alpha-\gamma}$ and $\eta_3:=(-1)^{\beta-\alpha}$, and in case $\varepsilon=-1$, $\eta_1:=(-1)^{\gamma-\beta}(-1)^{\chi_{\{k\leqslant p\}}+\chi_{\{j\leqslant p\}}}$, $\eta_2:=(-1)^{\alpha-\gamma}(-1)^{\chi_{\{i\leqslant p\}}+\chi_{\{k\leqslant p\}}}$, $\eta_3:=(-1)^{\beta-\alpha}(-1)^{\chi_{\{j\leqslant p\}}+\chi_{\{i\leqslant p\}}}$. We denote by $\chi_P$ the characteristic function of the set $P$, equal to 1 on $P$ and null elsewhere. Explicitly, $\chi_{\{i\leqslant p\}}+\chi_{\{j\leqslant p\}}$ is even if and only if $(i,j)\subset\llbracket1,p\rrbracket^2\cup \llbracket p+1, \delta\rrbracket^2$. Then the following sets of relations (say respectively {\bf (i)}, {\bf (ii)}, {\bf (iii)} and {\bf (iv)}) hold.
$$\left\{\begin{array}{l}
\eta_1\lambda_{H^{(1)}}(u',v,w)+\eta_2\lambda_{H^{(1)}}(u,v',w)\\\qquad\qquad\qquad
+\eta_3\lambda_{H^{(1)}}(u,v,w')+\varepsilon\lambda_{H^{(1)}}(u',v',w')=0\\[.9ex]
\eta_1\lambda_{H^{(1)}}(Ju',v,w)+\eta_2\lambda_{H^{(1)}}(u,Jv',w)\\\qquad\qquad\qquad
+\eta_3\lambda_{H^{(1)}}(u,v,Jw')+\varepsilon\lambda_{H^{(1)}}(Ju',Jv',Jw')=0\\[.9ex]
\eta_1\lambda_{H^{(1)}}(u,v',w')+\eta_2\lambda_{H^{(1)}}(u',v,w')\\\qquad\qquad\qquad
+\eta_3\lambda_{H^{(1)}}(u',v',w)+\varepsilon\lambda_{H^{(1)}}(u,v,w)=0\\[.9ex]
\eta_1\lambda_{H^{(1)}}(u,Jv',Jw')+\eta_2\lambda_{H^{(1)}}(Ju',v,Jw')\\\qquad\qquad\qquad
+\eta_3\lambda_{H^{(1)}}(Ju',Jv',w)+\varepsilon\lambda_{H^{(1)}}(u,v,w)=0.
\end{array}\right.$$
So the $8C_{2\delta}^3$ linear forms of the type $H^{(1)}\mapsto\lambda_{H^{(1)}}((J)u^{(\prime)},(J)v^{(\prime)},(J)w^{(\prime)})$ are linearly dependent, through the $4C_{2\delta}^3$ equations above. In turn, those equations are linearly independent. Counting the number of primes appearing in them, one sees that equations of types {\bf (i)}--{\bf (ii)} on the one hand, and types {\bf (iii)}--{\bf (iv)} on the other hand, span subspaces in direct sum. Now any dependence relation would involve some fixed triple $(i,j,k)$. For such a triple, equations of type {\bf (i)} may be seen as expressing the forms $H^{(1)}\mapsto\lambda_{H^{(1)}}((J)u_i^{\prime},(J)u_j^{\prime},(J)u_k^{\prime})$ as combination of the other ones, and then equations of type {\bf (i)}--{\bf (ii)}, doing the same with the forms $H^{(1)}\mapsto\lambda_{H^{(1)}}((J)u_i^{\prime},(J)u_j,(J)u_k)$. Equations of types {\bf (iii)}--{\bf (iv)} are similar, so all the $4C_{2\delta}^3$ equations are independent, and the $8C_{2\delta}^3$ forms span a space of dimension $\leqslant 8C_{2\delta}^3-4C_{2\delta}^3=4C_{2\delta}^3$. So {\bf (a)} and {\bf (b)} consist of not more than $4C_{2\delta}^3+4C_{2\delta}^2=4C_{2\delta+1}^3$ independent equations, so $\dim V_{4\delta}({\bf I})\geqslant\dim[\C^{2\delta}\otimes W_\varepsilon]-4C_{2\delta+1}^3=(4\delta).(2\delta^2+ \delta)-4C_{2\delta+1}^3=2C_{2\delta+2}^3$. This is {\bf (2)}.\medskip

We finally treat $\GG_\C$. In all that precedes, see all complex variables $x+{\rm i}y$ as real matrices $\text{\scriptsize$\left(\begin{array}{cc}x&y\\-y&x\end{array}\right)$}$. Then, complexifying everything {\em i.e.\@} replacing the real entries $x$, $y$ by complex numbers amounts to parametrise $\GG_\C$; so the same reasoning gives the proposition for $\GG_\C$.\hfill{\rm $\Box$}

\begin{importantrem}\label{generalite}Among real analytic germs of metrics with holonomy $H$ included in $H_0=\operatorname{Sp}(p,q)$, $H_0=\operatorname{Sp}(2\delta,\R)$ or $H_0=\operatorname{Sp}(2\delta,\C)$, corresponding to cases {\bf (3)}, {\bf (3')} and {\bf (3$^\C$)}, a dense open subset for the $C^2$ topology has its holonomy equal to $H_0$. Indeed, the first prolongation ${\bf I}^{(1)}$ of the ideal ${\bf I}$ satisfies also Cartan's criterion; this enables to show that any $2$-jet of metric, integrable at the order 1 and such that $\{R(X,Y);X,Y\in T_0\MM\}\subset {\goth h}_0$, is the 2-jet of a metric with holonomy included in $H_0$. The reasoning is presented, in the case $H=G_2$, in Proposition 3 p.\@ 556 of \cite{bryant1987}. It may be adapted here, as indicated in \cite{bryant1996} \S2.5 p.\@ 126. So as, among  such 2-jets, those satisfying  $\{R(X,Y);X,Y\in T_0\MM\}= {\goth h}_0$ are generic, we get the result.
\end{importantrem}

\section{Parallel endomorphisms and Ricci curvature}

The Ricci form $\ric(\,\cdot\,,J\,\cdot\,)$ has remarkable properties on K\"ahler manifolds. Let us determine the properties of the corresponding form $\ric(\,\cdot\,,U\,\cdot\,)$ when $g$ admits some other parallel endomorphism field $U$ than a K\"ahler structure.

\begin{te}\label{ricci} Suppose that $U$ is a parallel endomorphism field for a pseudo-Riemanian metric $g$; in the following $(a,b)$ denote any two tangent vectors at some point.\medskip

\noindent{\bf (i)} In case $U$ is self adjoint, the three following properties hold:
\begin{itemize}
\item[{\bf (a)}] $\ric(a,Ub)=\ric(Ua,b)=\tr(U(R(a,\cdot)b))$ and $U$ and $R(a,\cdot)b$ commute,

\item[{\bf (b)}] if $U=\underline J$ is a complex structure, then $g$ is the real part of the $\underline J$-complex metric $g_\C:=g(\cdot,\cdot)-{\rm i}g(\cdot,\underline J\cdot)$, and the Ricci curvature of $g_\C$ is $\ric_\C=\ric(\cdot,\cdot)-{\rm i}\ric(\cdot,\underline J\cdot)$,

\item[{\bf (c)}] if $U\!=\!N\!\neq0$ is nilpotent, $\ric$ is degenerate and $\im N\subset\ker\ric$.
\end{itemize}

\noindent{\bf (ii)} In case $U$ is skew adjoint, the three following properties hold:
\begin{itemize}
\item[{\bf (a)}] $\ric(a,Ub)=-\ric(Ua,b)=\frac12\tr(U\circ R(a,b))$,

\item[{\bf (b)}] if $U=N\neq0$ is nilpotent, $\ric$ is degenerate and $\im N\subset\ker\ric$,

\item[{\bf (c)}] if $V$ is another skew symmetric parallel endomorphism, $\im (UV-VU)\subset\ker\ric$. In particular, if moreover $U$ and $V$ are invertible and anti-commute, $\ric=0$. It follows that cases {\bf (3)}, {\bf (3')}, {\bf (3$^\C$)} of Theorem \ref{structure_s} are Ricci-flat.
\end{itemize}
\end{te}

\begin{rem}
Point {\bf (i)(b)} is a standard result and  we will not prove it. The Ricci flatness of cases {\bf (3)}, {\bf (3')}, {\bf (3$^\C$)} is also standard and classically proven by other means, we show here that it can be deduced elementarily from {\bf (ii)(a)}.
\end{rem}

\begin{rem}
The space $\goth e^-$ of parallel skew adjoint endomorphism fields is the Lie algebra $\goth o(g)\cap\goth e$. Point {\bf (ii)(c)} means that $[\goth e^-,\goth e^-]\subset\ker\ric$.
\end{rem}

\noindent{\bf Proof.} Take $U$ self adjoint, then the whole of {\bf (a)} follows from Remark \ref{g_U}. For {\bf (c)}, after {\bf (a)}, $\ric(a,Nb)=\tr(N(R(a,\cdot)b))$, and as $N$ and $R(a,\,\cdot\,)b$ commute, their product is also nilpotent, so trace free. Now take $U$ skew adjoint.
\begin{align*}
\ric(a,Ub)&=\tr(R(a,\,\cdot\,)Ub)\\
&=\tr(U(R(a,\,\cdot\,)b))\ \text{as $U$, being parallel, commutes with $R(a,\,\cdot\,)$,}\\
&=\tr(R(a,U\,\cdot\,)b)\quad\text{as $\tr(UV)=\tr(VU)$,}\\
&=-\tr(R(Ua,\,\cdot\,)b).
\end{align*}
For the last line, take any $u$, $v$, $w$: $g(R(Ua,u)v,w)=g(R(v,w)Ua,u)=g(UR(v,w)a,u)=-g(R(v,w)a,Uu)=-g(R(a,Uu)v,w)$. So finally, $\ric(a,Ub)=-\ric(Ua,b)$. Besides:
\begin{align*}
\ric(Ua,b)&=\ric(b,Ua)&&\text{by the first equality of {\bf (a)},}\\
&=\tr(U(R(b,\,\cdot\,)a))&&\text{by definition of $\ric$,}\\
&=-\tr(U(R(\,\cdot\,,a)b))-\tr(U(R(a,b)\,\cdot\,))
&&\text{by the Bianchi identity,}\\
&=\tr(R(a,\,\cdot\,)Ub)-\tr(U\circ R(a,b))&&\text{as $U$ commutes with $R(a,\,\cdot\,)=-R(\,\cdot\,,a)$,}\\&=\ric(a,Ub)-\tr(U\circ R(a,b)).
\end{align*}
As $\ric(Ua,b)=-\ric(a,Ub)$, we get {\bf (a)}. Point {\bf (b)} follows: $\ric(a,Nb)=\frac12\tr(N\circ R(a,b))=0$ as $N\circ R(a,b)=R(a,b)\circ N$ is nilpotent. We are left with proving the first assertion of {\bf (c)} {\em i.e.\@} that $\ric(a,UVb)=\ric(a,VUb)$.
\begin{align*}
\ric(a,UVb)&=-\ric(Ua,Vb)&&\text{by {\bf (a)} applied to $U$,}\\
&=-{\textstyle\frac12}\tr(V\circ R(Ua,b))&&\text{by {\bf (a)} applied to $V$,}\\
&={\textstyle\frac12}\tr(V\circ R(a,Ub))&&\text{as $U^\ast=-U$,}\\
&=\ric(a,VUb)&&\text{by {\bf (a)} applied to $V$.}\tag*{$\Box$}
\end{align*}

\begin{cor}$\!$Let $\Ric$ be the endomorphism such that $\ric=g(\,\cdot\,,\Ric\,\cdot\,)$. If the metric is indecomposable (in a local Riemannian product) and such that $\ric$ is parallel, then $\Ric$ is either semi-simple or 2-step nilpotent.
\end{cor}

\noindent{\bf Proof.} As $g$ is indecomposable, the minimal polynomial of $\Ric$ is of the form $P^\alpha$ with $P$ irreductible, see {\em Claim 1} p.\@ \pageref{polynome_minimal} in the proof of Theorem \ref{structure_s}. So $\Ric$ is either invertible or nilpotent. Apply Theorem \ref{ricci} {\bf (i)} {\bf (c)} to the nilpotent part $N_{\Ric}$ of $\Ric$: if $\Ric$ is invertible, $\ker\ric=\{0\}$ so $N_{\Ric}=0$, else $\Ric^2=N_{\Ric}^2=0$. We re-find here the result of \cite{boubel-berard}.\hfill{\rm $\Box$}


{\footnotesize
\begin{thebibliography}{0}
\bibitem{BB-ikemakhen1997} {\sc L.\@ B\'erard-Bergery, A.\@ Ikemakhen},  Sur l'holonomie des va\-ri\'e\-t\'es pseudo-riemanniennes de signature $(n,n)$. {\em Bull.\@ Soc.\@ Math.\@ France} {\bf 125} no. 1, 93--114 (1997).
\bibitem{boubel13}{\sc C.\@ Boubel}, The algebra of the parallel endomorphisms of a germ of pseudo-Riemannian metric. {\em J.\@ Differential Geom.\@} {\bf 99} No.\@ 1, 77--123 (2015).
\bibitem{boubel-berard}{\sc C.\@ Boubel, L.\@ B\'erard-Bergery} On pseudo-Riemannian manifolds whose Ricci tensor is parallel. {\em Geom.\@ Dedicata} {\bf 86}, No.1-3, 1--18 (2001).
\bibitem{bourbaki_algebre8}{\sc N.\@ Bourbaki}, {\em \'El\'ements de math\'ematique. Livre II: Alg\`ebre. Chap.\@ 8: Modules et anneaux semisimples.} (French) Actualit\'es scientifiques et industrielles, Hermann, 1958.
\bibitem{bryant1987} {\sc R.\@ Bryant}, Metrics with exceptional holonomy. {\em Ann.\@ Math.\@} {\bf 126}, 525--576 (1987).
\bibitem{bcggg} {\sc R.\@  L.\@ Bryant, S.\@  S.\@ Chern, R.\@  B.\@ Gardner, H.\@  L.\@ Goldschmidt, P.\@  A.\@ Griffiths}, {\em  Exterior differential systems.}  Mathematical Sciences Research Institute Publications, 18.  Springer-Verlag, 1991.
\bibitem{bryant1996} {\sc R.\@ Bryant}, Classical, exceptional, and exotic holonomies: A status report. In: {\em A.\@ Besse (ed.), Actes de la table ronde de g\'eom\'etrie diff\'erentielle en l'honneur de Marcel Berger, Luminy, France, 12--18 juillet, 1992.}, 93--165. Soci\'et\'e Math\'ematique de France, S\'emin.\@ Congr., 1, 1996.
\bibitem{bryant2001} {\sc R.\@ Bryant}, Bochner-K\"ahler metrics, {\em J.\@ Amer.\@ Math.\@ Soc.\@} {\bf 14} No.\@ 3, 623-715 (2001).
\bibitem{curtis-reiner}{\sc C.\@ Curtis, I.\@ Reiner}, {\em Representation theory of finite groups and associative algebras. Reprint of the 1962 original.} AMS Chelsea Publishing, 2006.
\bibitem{galaev-leistner2010} {\sc A.\@ Galaev, T.\@ Leistner}, Recent developments in pseudo-Riemannian holonomy theory. In: {\em Handbook of pseudo-Riemannian geometry and supersymmetry}, 581--627. IRMA Lect.\@ Math.\@ Theor.\@ Phys.\@, 16, Eur.\@ Math.\@ Soc.\@, 2010.
\bibitem{ghanam-thompson2001} {\sc R.\@ Ghanam, G.\@ Thompson}, The holonomy Lie algebras of neutral metrics in dimension four. {\em J.\@ Math.\@ Phys.\@} {\bf 42} no. 5, 2266--2284 (2001).
\bibitem{ivey-landsberg} {\sc T.\@ A.\@ Ivey, J.\@ M.\@ Landsberg}, {\em Cartan for beginners: differential geometry via moving frames and exterior differential systems.} Graduate Studies in Mathematics, 61. Amer.\@ Math.\@ Soc.\@, 2003.
\bibitem{jacobson_structureofrings}{\sc N.\@ Jacobson}, {\em Structure of rings.}, AMS, Colloquium Publications, vol.\@ 37, 1956.
\bibitem{milnor_morse}{\sc J.\@ Milnor}, {\em Morse theory. Based on lecture notes by M. Spivak and R. Wells.} Annals of Mathematics Studies no. 51. Princeton University Press,1963.
\bibitem{moroianu}{\sc A.\@ Moroianu}, {\em Lectures on K\"ahler Geometry}, Cambridge University Press 2007.

\bibitem{schwachhoefer} {\sc L.\@ Schwachh\"ofer}, Connections with irreducible holonomy representations, {\em Adv.\@ Math.\@} {\bf 160} No. 1, 1--80 (2001).
\bibitem{taft1} {\sc E.\@ J.\@ Taft}, Invariant Wedderburn factors, {\em Illinois J.\@ Math.\@} {\bf 1}, 565--573 (1957).
\bibitem{taft2} {\sc E.\@ J.\@ Taft}, Cleft algebras with operator groups, {\em Portugal.\@ Math.\@} {\bf 20}, 195--198 (1961).
\end{thebibliography}}
\end{document}